\documentclass[11pt]{amsart}
\usepackage[letterpaper,margin=1in]{geometry}               
\usepackage{amssymb,amsmath,amsfonts,stmaryrd}
\usepackage{xcolor}
\usepackage[colorlinks,
             linkcolor=blue,
             citecolor=green!70!black,
             pdfproducer={pdfLaTeX},
             pdfpagemode=None,
             bookmarksopen=true
             bookmarksnumbered=true]{hyperref}
\usepackage{graphicx}
\usepackage{placeins}
\usepackage{multicol}
\usepackage{tikz-cd}
\usetikzlibrary{shapes.geometric,calc}
\usetikzlibrary{shapes,arrows,chains}
\usetikzlibrary{decorations.markings}
\usetikzlibrary{decorations.pathmorphing}
\usetikzlibrary{shapes.multipart}
\tikzset{snake it/.style={decorate, decoration=snake}}
\tikzstyle{GraphNode}=[circle, draw=black, fill=black, inner sep=2pt, minimum size=5pt]
\tikzstyle{GraphEdge}=[black]
\usepackage{lipsum}
\usepackage{url}
\usepackage{hyperref}
\usepackage{adjustbox}
\usepackage{float}
\usepackage{tikz}
\usepackage{bbold}
\usepackage{tikz-cd}
\usetikzlibrary{shapes.geometric,decorations.markings,arrows,decorations.pathreplacing,calc}
\usepackage{arydshln}
\newcommand{\boundellipse}[3]
{(#1) ellipse (#2 and #3)
}

\usepackage{comment}

\newcommand\arXiv[1]{\href{http://arxiv.org/abs/#1}{\nolinkurl{arXiv:#1}}}
\newcommand\MRnumber[1]{\href{http://www.ams.org/mathscinet-getitem?mr=#1}{\nolinkurl{MR#1}}}
\newcommand\DOI[1]{\href{http://dx.doi.org/#1}{\nolinkurl{DOI:#1}}}
\newcommand\MAILTO[1]{\href{mailto:#1}{\nolinkurl{#1}}}

\newcounter{mainthm}

\newtheorem{maintheorem}[mainthm]{Theorem}
\newtheorem{dummy}{Dummy}[section]

\newtheorem{lemma}[dummy]{Lemma}
\newtheorem{proposition}[dummy]{Proposition}
\newtheorem{corollary}[dummy]{Corollary}
\newtheorem{definition}[dummy]{Definition}
\newtheorem{theorem}[dummy]{Theorem}

\theoremstyle{definition}

\newtheorem*{rem}{Remark}

\newtheorem*{example}{Example}

\usepackage{dsfont}
\renewcommand\mathbb\mathds

\newcommand\bC{\mathbb C}

\newcommand\bZ{\mathbb Z}

\newcommand\rH{\mathrm H}

\newcommand\SH{\mathrm {SH}}

\newcommand\SW{\mathrm{S}\mathcal{W}}

\newcommand\pt{\mathrm{pt}}
\usepackage[mathscr]{euscript}

\newcommand\fB{\mathfrak B}
\newcommand\fC{\mathfrak C}
\newcommand\fS{\mathfrak S}

\newcommand\longto\longrightarrow
\newcommand\mono\hookrightarrow
\newcommand\epi\twoheadrightarrow

\newcommand\<\langle
\renewcommand\>\rangle
\newcommand\sminus\smallsetminus

\DeclareMathOperator\End{End}
\DeclareMathOperator{\Sq}{Sq}
\DeclareMathOperator\Ext{Ext}
\DeclareMathOperator\Quad{Quad}

\DeclareMathOperator{\Spin}{Spin}

\DeclareMathOperator{\Gal}{Gal}

\newcommand\SVec{\cat{SVec}}
\renewcommand\Vec{\cat{Vec}}
\newcommand\Mod{\cat{Mod}}
\newcommand\rRep{\cat{Rep}}

\newcommand\define[1]{\emph{#1}}
\newcommand\cat[1]{\mathbf{#1}}
\newcommand\matt[1]{{\textcolor{purple}{*#1*}}}
\newcommand\thib[1]{{\textcolor{red}{*#1*}}}

\title{Gauging Noninvertible Defects: A 2-Categorical Perspective} 
\author[D\'ecoppet and Yu]{Thibault D. D\'ecoppet$^{\lowercase{a}}$ and Matthew Yu$^{\lowercase{b}}$}
\thanks{It is a pleasure to thank Arun Debray, Jonathan Heckman, and Theo Johnson-Freyd for helpful conversations and comments. We would also like to thank the Simons Collaboration on Global Categorical Symmetries for hosting a summer school where this work began. Research at the Perimeter Institute is
supported by the Government of Canada through Industry Canada and by the Province of Ontario through
the Ministry of Economic Development and Innovation.
 \\[6pt]
$^a$ \textsc{Mathematical Institute, University of Oxford, United Kingdom}.\\ 
$^b$ \textsc{Perimeter Institute for Theoretical Physics, Waterloo, Ontario}.
\\[6pt]
}

\begin{document}

\maketitle

\begin{abstract}
    We generalize the notion of an anomaly for a symmetry to a noninvertible symmetry enacted by surface operators using the framework of condensation in 2-categories. Given a multifusion 2-category, potentially with some additional levels of monoidality, we prove theorems about the structure of the 2-category obtained by condensing a suitable algebra object. We give examples where the resulting category displays grouplike fusion rules and through a cohomology computation, find the obstruction to condensing further to the vacuum theory. As a consequence, we show that every symmetric fusion 2-category admits a fibre 2-functor to $2\SVec$.
\end{abstract}

\section{Introduction}

One of the most exciting prospects of generalized symmetries is the study of noninvertible symmetry operators. These are topological but instead of having a grouplike composition, their interactions are described by a general higher category. For grouplike global symmetries, the anomaly determines whether or not the symmetry can be gauged. The classification of such anomalies is well known to be captured by an invertible theory one dimension higher. Further, they can be classified using spectral sequences for group cohomology, and, more generally, for cobordisms as formalized in \cite{Freed:2016,Kapustin:2014}. Provided the anomaly vanishes, the gauging procedure will in general reshuffle the topological content, and in some cases add new richness into the theory in form of noninvertible operators \cite{Kaidi:2021,Choi:2021}. 
When gauging discrete abelian groups, what manifests is a dual group, which upon gauging takes us back to the original theory. The notion of \textit{condensation} was introduced in \cite{GJF} as a generalization of gauging, which applies to noninvertible symmetries. One particularly useful perspective of condensing a symmetry involves starting from the vacuum theory and proliferating in space (or perhaps in some subspace) a network of operators for that symmetry which fill out a new phase \cite{Roumpedakis:2022aik,Choi:2022zal}. Since this procedure is fully topological one can imagine running this procedure backwards and constructing a topological boundary between some phase and the vacuum. If one can go back and forth with no obstruction, then the symmetry is nonanomalous.

The purpose of this article is to generalize the notion of an anomaly for a symmetry, to an anomaly for a noninvertible symmetry. We will focus on noninvertible surface operators, for which the natural mathematical setting is a 2-category. For other applications of 2-categories in the physics literature, we refer the reader to \cite{Bhardwaj:2022yxj, Bartsch:2022, Bhardwaj:universal}. In general, the 2-category $\mathfrak{C}$ can have more structure such as a braiding, where the braiding takes place along the morphisms of $\mathfrak{C}$, or a syllepsis, and we will consider both cases.
If one is in a setting were the surfaces are fully symmetric, we will show that a higher analogue of Deligne's theorem in \cite{deligne2002} holds. More precisely, it was first announced in \cite{JF3}, that for any symmetric fusion 2-category $\fS$, there exists a fibre 2-functor $\text{Fib}:\fS \to 2\SVec$ to the 2-category of super-2-vector spaces. In this sense, in the fully symmetric case, there is no obstruction to condensing all the operators, if we allow for emergent fermions.
In this work we will consider both the cases of condensing to $2\Vec$, the 2-category of 2-vector spaces, and to $2\SVec$, where the latter involves working fermionically by condensing a fermionic algebra. This is the noninvertible analogue of being able to gauge a symmetry.
In this article, we are mainly concerned with theories that have surface operators belonging to a fusion 2-category $\mathfrak{C}$ that can at least braid with each other, but are not fully symmetric. Since $\fC$ is not fully symmetric, there is no universal target that all the operators can condense to. 
We instead consider a related question which involves finding a subcategory of surface operators that enjoy more levels of monoidality than the general surface operators in the ambient category. One such example is given by the
extra data of the aforementioned syllepsis, which can be thought of as anomaly cancellation data associated to the braiding \footnote{The additional level of monoidality means that the surface can secretly ascend to a higher dimension. For example, surfaces can braid in four total dimensions, but the data of being sylleptic means that some set of surfaces can lift to five dimensions.}.
It is then a meaningful question to ask what happens to the ambient category upon condensing the subcategory. While working in a 2-category, if it so happens that there exists a procedure to go to the vacuum theory, then there will be no anomalies for any noninvertible symmetry, as all of them will have been ``gauged". This idea will be useful in theories of gravity where it is expected that, not only there are no global symmetries, but also no noninvertible symmetries. For more on global symmetries arising in gravitational settings see \cite{Apruzzi:2022,Heckman:2022,GarciaEtxebarria:2022,Albertini:2020,vanBeest:2022}.

Building on the work of the first author \cite{D7, D8}, the main results of this article are proven in \S\ref{section:BraidedSymmetric}. More precisely, we present the result of condensing noninvertible surfaces in an ambient 2-category, with subsequent corollaries involving changing the properties of the condensation monad, also called separable algebra.

\begin{maintheorem}\label{thm:first}
For $\mathfrak{B}$ a braided multifusion 2-category, condensing a braided separable algebra $B$ in $\mathfrak{B}$ results in a multifusion 2-category.
\end{maintheorem}

\begin{maintheorem}\label{thm:second}
For $\mathfrak{S}$ a sylleptic multifusion 2-category, condensing a symmetric separable algebra $B$ in $\mathfrak{S}$, results in a braided multifusion 2-category.
\end{maintheorem}

\begin{maintheorem}\label{thm:third}
For $\mathfrak{S}$ a sylleptic multifusion 2-category, condensing a symmetric separable algebra $B$ in the symmetric center of $\mathfrak{S}$, results in a sylleptic multifusion 2-category. Further, if $\mathfrak{S}$ is symmetric, then condensing $B$ yields a symmetric multifusion 2-category.
\end{maintheorem}

\noindent The auxiliary results of this article build off the main theorems by exploring particularly nice cases where the resulting category after condensation is ``grouplike", in addition to being braided, or sylleptic. We call these categories \textit{strongly fusion}, and the operator content is essentially captured by the surfaces \cite{JFY}.
The reader interested in applications of the main theorems can go to \S\ref{section:ExampleCats} for explicit examples of condensations within 2-categories, which in the right setting, yield strongly fusion categories. In particular, we show that every symmetric fusion 2-category can be condensed to a symmetric strongly fusion 2-category.
For theories described by strongly fusion 2-categories, the obstruction to condensing to the vacuum is given by a cohomology class, which we compute when the 2-category is braided. In addition, we show that the obstruction to condense a symmetric strongly fusion 2-category to the 2-category of super-2-vector spaces vanishes. Thereby establishing the following result:

\begin{maintheorem}\label{main:fourth}
Every symmetric fusion 2-category admits a fibre 2-functor to $2\SVec$.
\end{maintheorem}

\noindent The above theorem categorifies \cite{deligne2002}.

We now outline the contents of this article:
In \S\ref{section:preliminary} we explain the graphical calculus used for braided, sylleptic, and symmetric monoidal 2-categories. We also discuss algebras, and the relationship between modules and condensation. 
In \S\ref{section:BraidedSymmetric} we prove the main theorems about braided or sylleptic monoidal 2-category, and the result of condensing separable algebras that are respectively braided or symmetric. We examine specific examples of condensing separable algebras in connected and disconnected 2-categories that are interesting for physical applications in \S\ref{section:ExampleCats}; we find that in some cases, the 2-category becomes \textit{strongly fusion}. In \S\ref{section:StronglyFusion} we perform cohomology computations for theories described by 
the braided and symmetric strongly fusion 2-categories, and report on the obstruction to condensing the theory to the vacuum.

\section{Preliminaries on 2-Categories}\label{section:preliminary}
\subsection{Graphical Calculus} We begin by setting up the fundamental definitions and explaining the computational language of string diagrams. We work within a monoidal 2-category $\mathfrak{C}$ with monoidal unit $I$ and monoidal product $\Box$ in the sense of definition 2.3 of \cite{SP}. Thanks to the coherence theorem of \cite{Gur}, we may assume without loss of generality that $\mathfrak{C}$ is strict cubical (in the sense that it satifiesthe conditions of definition 2.26 of \cite{SP}). In this setting, we use the graphical calculus of \cite{GS}, as described in \cite{D7} (see also \cite{D4}). In particular, we will often omit the monoidal product $\Box$ from our notation. In addition, identity 1-morphisms are denoted using the symbol $1$. Further, the interchanger is depicted using by the string diagram below on the left, and its inverse by that on the right: $$\begin{tabular}{c c c c}
\includegraphics[width=20mm]{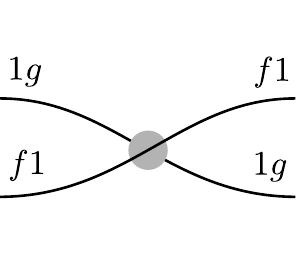},\ \ \ \   & \ \ \ \  \includegraphics[width=20mm]{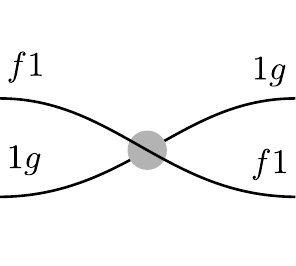}.
\end{tabular}$$
The lines represent 1-morphisms and their composition is read from top to bottom. The string diagrams are then read from left to right, and the coupons represent 2-morphisms.
The regions between the lines represent objects of the 2-category, which are specified uniquely by the 1-morphisms. 

We also need to recall the graphical conventions related to 2-natural transformations from \cite{GS}. In the present article, these will exclusively be used for the braiding, which will be introduced below. Let $F,G:\mathfrak{A}\rightarrow \mathfrak{B}$ be two (weak) 2-functors, and let $\tau:F\Rightarrow G$ be 2-natural transformation. This means that, for every object $A$ in $\mathfrak{A}$, we have a 1-morphism $\tau_A:F(A)\rightarrow G(A)$, and for every 1-morphism $f:A\rightarrow B$ in $\mathfrak{A}$, we have a 2-isomorphism $$\begin{tikzcd}[sep=tiny]
F(A) \arrow[ddd, "F(f)"']\arrow[rrr, "\tau_A"]  &                                        &    & G(A) \arrow[ddd, "G(f)"]  \\
 &  &    & \\
  &  &  &  \\
F(B)\arrow[rrr, "\tau_B"']\arrow[Rightarrow, rrruuu, "\tau_f", shorten > = 2ex, shorten < = 2ex]                                            &                                        &    &  G(B), 
\end{tikzcd}$$
The collection of these 2-isomorphisms has to satisfy the obvious coherence relations. In our graphical language, we will depict the 2-isomorphism $\tau_f$ using the following diagram on the left, and its inverse using the diagram on the right: $$\begin{tabular}{c c c c}
\includegraphics[width=20mm]{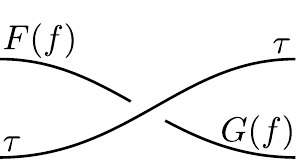},\ \ \ \   & \ \ \ \  \includegraphics[width=20mm]{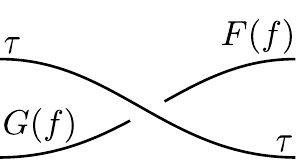}.
\end{tabular}$$ 

\subsubsection{Braided Monoidal 2-Categories}\label{subsub:braidedmonoidal2cat}

Let $\mathfrak{B}$ a braided monoidal 2-category in the sense of definition 2.3 of \cite{SP}. In particular, $\mathfrak{C}$ is a monoidal 2-category, so that we use $I$ to denote its monoidal unit, and $\Box$ to denote its monoidal product. The coherence theorem of \cite{Gur2} allows us to assume that $\mathfrak{B}$ is a semi-strict braided monoidal 2-category. In particular, the underlying monoidal 2-category is strict cubical. Further, $\mathfrak{B}$ comes equipped with a braiding $b$, which is an adjoint 2-natural equivalence given on objects $A,B$ in $\mathfrak{B}$ by $$b_{A,B}:A\Box B\rightarrow B\Box A.$$ Further, there are two invertible modifications $R$ and $S$, which are given on the objects $A,B,C$ of $\mathfrak{B}$ by \begin{center}
\begin{tabular}{@{}c c@{}}

$\begin{tikzcd}
ABC \arrow[rr, "b"] \arrow[rd, "b1"'] & {} \arrow[d, Rightarrow, "R"]          & BCA, \\
                                      & BAC \arrow[ru, "1b"'] &    
\end{tikzcd}$

&

$\begin{tikzcd}
ABC \arrow[rr, "b_2"] \arrow[rd, "1b"'] & {} \arrow[d, Rightarrow, "S"]     & CAB \\
                                        & ACB \arrow[ru, "b1"'] &    
\end{tikzcd}$
\end{tabular}
\end{center}
where the subscript in $b_2$ records that the braiding occurs between the first two objects on the left and the next ones. On the other hand, $b$ means that the braiding occurs between the first object on the left and the next ones. These two modifications are subject to the following relations:

\begin{enumerate}
\item [a.] For all objects $A,B,C,D$ in $\mathfrak{B}$, we have
\end{enumerate}

\newlength{\calculus}

\settoheight{\calculus}{\includegraphics[width=37.5mm]{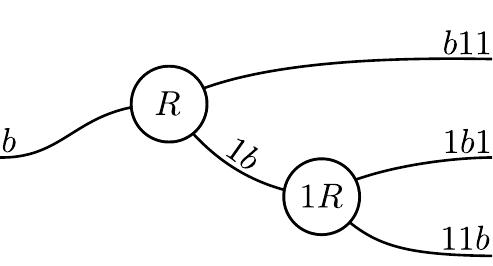}}

\begin{equation}\label{eqn:braidingaxiom1}
\begin{tabular}{@{}ccc@{}}

\includegraphics[width=37.5mm]{Pictures/Preliminaries/Braiding/braidingaxiom1.pdf} & \raisebox{0.45\calculus}{$=$} &
\includegraphics[width=37.5mm]{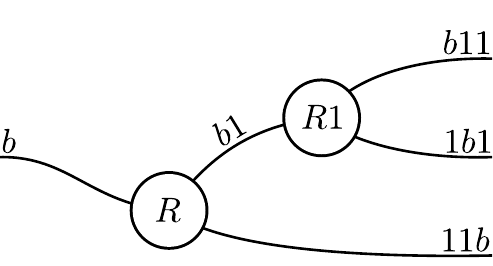}

\end{tabular}
\end{equation}
\begin{enumerate}
\item [] in $Hom_{\mathfrak{B}}(ABCD, BCDA)$. 
\item [b.] For all objects $A,B,C,D$ in $\mathfrak{B}$, we have
\end{enumerate}

\settoheight{\calculus}{\includegraphics[width=37.5mm]{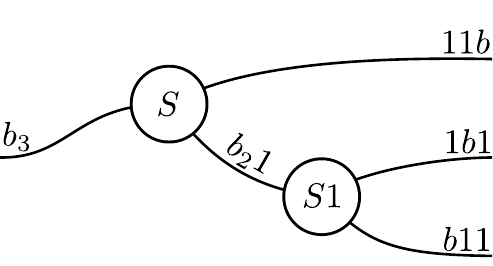}}

\begin{equation}\label{eqn:braidingaxiom2}
\begin{tabular}{@{}ccc@{}}

\includegraphics[width=37.5mm]{Pictures/Preliminaries/Braiding/braidingaxiom3.pdf} & \raisebox{0.45\calculus}{$=$} &
\includegraphics[width=37.5mm]{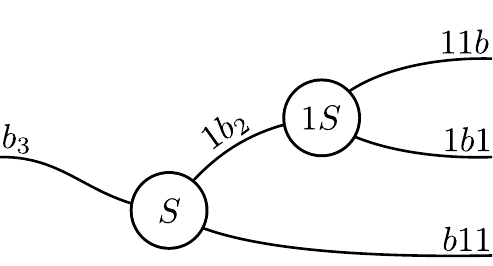}

\end{tabular}
\end{equation}

\begin{enumerate}
\item [] in $Hom_{\mathfrak{B}}(ABCD, DABC)$,
\item [c.] For all objects $A,B,C,D$ in $\mathfrak{B}$, we have
\end{enumerate}

\settoheight{\calculus}{\includegraphics[width=37.5mm]{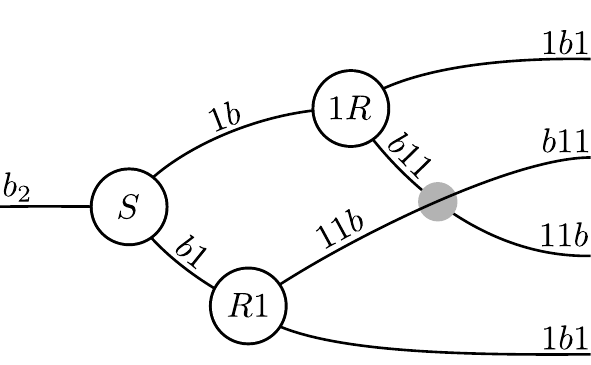}}

\begin{equation}\label{eqn:braidingaxiom3}
\begin{tabular}{@{}ccc@{}}

\includegraphics[width=45mm]{Pictures/Preliminaries/Braiding/braidingaxiom5.pdf} & \raisebox{0.45\calculus}{$=$} &
\includegraphics[width=37.5mm]{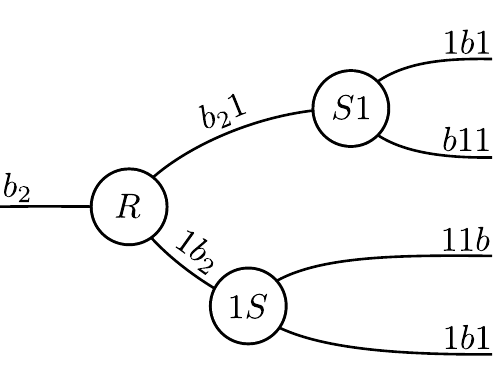}

\end{tabular}
\end{equation}

\begin{enumerate}
\item [] in $Hom_{\mathfrak{B}}(ABCD, CDAB)$,
\item [d.] For all objects $A,B,C$ in $\mathfrak{B}$, we have
\end{enumerate}

\settoheight{\calculus}{\includegraphics[width=45mm]{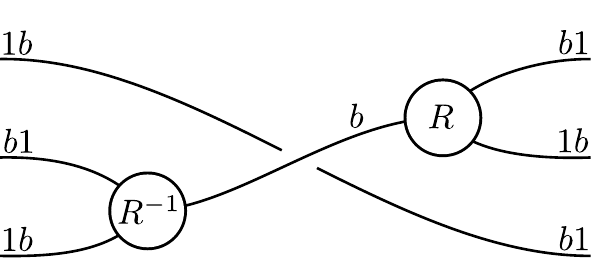}}

\begin{equation}\label{eqn:braidingaxiom4}
\begin{tabular}{@{}ccc@{}}

\includegraphics[width=45mm]{Pictures/Preliminaries/Braiding/braidingaxiom7.pdf} & \raisebox{0.45\calculus}{$=$} &
\includegraphics[width=45mm]{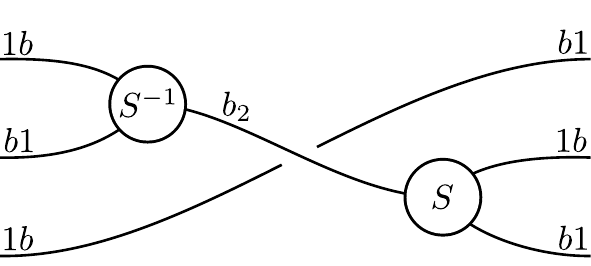}

\end{tabular}
\end{equation}

\begin{enumerate}
\item [] in $Hom_{\mathfrak{B}}(ABC, CBA)$,
\item [e.] For all objects $A$ in $\mathfrak{B}$, the adjoint 2-natural equivalences $$b_{A,I}:A\Box I\rightarrow I\Box A\textrm{ and } b_{I,A}:I\Box A\rightarrow A\Box I$$ are the identity adjoint 2-natural equivalences,

\item [f.] For all objects $A,B,C$ in $\mathfrak{B}$, the 2-isomorphisms $R_{A,B,C}$ and $S_{A,B,C}$ are equal to the identity 2-isomorphism whenever either $A$, $B$, or $C$ is equal to $I$.
\end{enumerate}
 In each of the $Hom_{\mathfrak{B}}$ above, the first set of objects is given by the top most region bound by 1-morphism, and the second set of objects is given by the bottom most region.

\subsubsection{Sylleptic and Symmetric Monoidal 2-Categories} \label{subsub:syllepticmonoidal2cat}
Our work will also involve sylleptic monoidal 2-categories (see definition 2.3 of \cite{SP}), these are braided monoidal 2-categories equipped with an additional structure called a \textit{syllepsis}. Without loss of generality, we may assume that every sylleptic monoidal 2-category $\mathfrak{S}$ semi-strict. (This follows from a slight generalization of \cite{GO}.) This means that $\mathfrak{S}$ is a semi-strict braided monoidal 2-category equipped with an invertible modification $\sigma$ given on the objects $A,B$ of $\mathfrak{S}$ by  $$\begin{tikzcd}[sep=small]
AB \arrow[rrrr, equal] \arrow[rrdd, "b"'] &  & {} \arrow[dd, Rightarrow, "\sigma", near start, shorten > = 1ex] &  & AB. \\
                                   &  &                           &  &   \\
                                   &  & BA \arrow[rruu, "b"']     &  &  
\end{tikzcd}$$

\noindent Furthermore, the invertible modification $\sigma$ satisfies the following relations:

\begin{enumerate}
\item [a.] For all objects $A,B,C$ of $\mathfrak{S}$, we have
\end{enumerate}

\settoheight{\calculus}{\includegraphics[width=30mm]{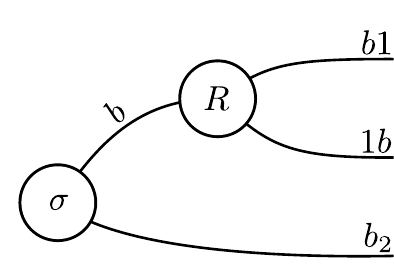}}

\begin{equation}\label{eqn:syllepsisaxiom1}
\begin{tabular}{@{}ccc@{}}

\includegraphics[width=30mm]{Pictures/Preliminaries/Braiding/syllepsisaxiom1.pdf} & \raisebox{0.45\calculus}{$=$} &
\includegraphics[width=45mm]{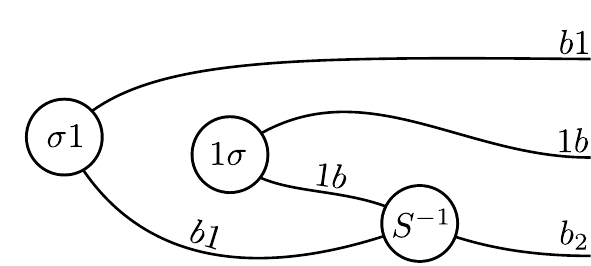}

\end{tabular}
\end{equation}

\begin{enumerate}
\item [] in $Hom_{\mathfrak{B}}(ABC, ABC)$,
\item [b.] For all objects $A,B,C$ of $\mathfrak{S}$, we have
\end{enumerate}

\settoheight{\calculus}{\includegraphics[width=30mm]{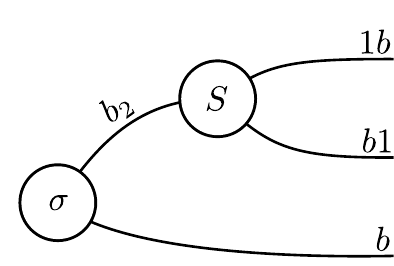}}

\begin{equation}\label{eqn:syllepsisaxiom2}
\begin{tabular}{@{}ccc@{}}

\includegraphics[width=30mm]{Pictures/Preliminaries/Braiding/syllepsisaxiom3.pdf} & \raisebox{0.45\calculus}{$=$} &
\includegraphics[width=45mm]{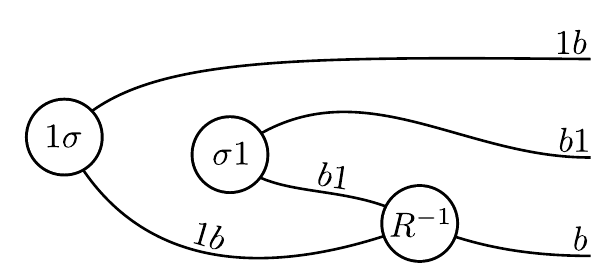}
\end{tabular}
\end{equation}

\begin{enumerate}\label{eq:conditionC}
\item [] in $Hom_{\mathfrak{B}}(ABC, ABC)$,
\item [c.] For all objects $A,B$ of $\mathfrak{S}$, the 2-isomorphisms $\sigma_{A,B}$ is the identity 2-morphism whenever either $A$ or $B$ is equal to $I$.
\end{enumerate}
We give a physical interpretation of syllepsis for surfaces. Namely, two surfaces existing in 5d braid by passing one around around each other in a 2 parameter family. The surfaces can exchange the order of which one is on top by going into the fifth dimension and using the syllepsis.

Finally, we will also consider symmetric monoidal 2-categories. Thanks to the main result of \cite{GO}, every symmetric monoidal 2-category is equivalent to a semi-strict symmetric monoidal 2-category that is to a semi-strict sylleptic monoidal 2-category $\mathfrak{S}$ as defined above, whose syllepsis satisfies \begin{equation}\label{eqn:symmetricaxiom}\sigma_{B,A}\circ b_{A,B} = b_{A,B}\circ \sigma_{A,B}\,,\end{equation} for every object $A,B$ in $\mathfrak{S}$. 
Physically speaking, if the surface operators have enough freedom to move around each other, such as in six ambient spacetime dimensions, then this is automatic. 

\subsection{Algebras and Modules}

Let $\mathfrak{C}$ be a strict cubical monoidal 2-category. We recall the definition of an algebra in $\mathfrak{C}$ expressed using our graphical calculus from \cite{D7}. These objects were introduced under the name pseudo-monoidal in \cite{DS}. The definition of an algebra in an arbitrary monoidal 2-category using our graphical conventions may be found in \cite{D4}.
\begin{definition}\label{def:algebra}
An algebra in $\mathfrak{C}$ consists of:
\begin{enumerate}
    \item An object $A$ of $\mathfrak{C}$;
    \item Two 1-morphisms $m:A\Box A\rightarrow A$ and $i:I\rightarrow A$;
    \item Three 2-isomorphisms
\end{enumerate}
\begin{center}
\begin{tabular}{@{}c c c@{}}
$\begin{tikzcd}[sep=small]
A \arrow[rrrr, equal] \arrow[rrdd, "i1"'] &  & {} \arrow[dd, Rightarrow, "\lambda"', near start, shorten > = 1ex] &  & A \\
                                   &  &                           &  &   \\
                                   &  & AA, \arrow[rruu, "m"']     &  &  
\end{tikzcd}$

&

$\begin{tikzcd}[sep=small]
AAA \arrow[dd, "1m"'] \arrow[rr, "m1"]    &  & AA \arrow[dd, "m"] \\
                                            &  &                      \\
AA \arrow[rr, "m"'] \arrow[rruu, Rightarrow, "\mu", shorten > = 2.5ex, shorten < = 2.5ex] &  & A,                   
\end{tikzcd}$

&

$\begin{tikzcd}[sep=small]
                                  &  & AA \arrow[rrdd, "m"] \arrow[dd, Rightarrow, "\rho", shorten > = 1ex, shorten < = 2ex] &  &   \\
                                  &  &                                             &  &   \\
A \arrow[rruu, "1i"] \arrow[rrrr,equal] &  & {}                                          &  & A,
\end{tikzcd}$

\end{tabular}
\end{center}

satisfying:

\begin{enumerate}
\item [a.] 
\end{enumerate}

\newlength{\prelim}

\settoheight{\prelim}{\includegraphics[width=52.5mm]{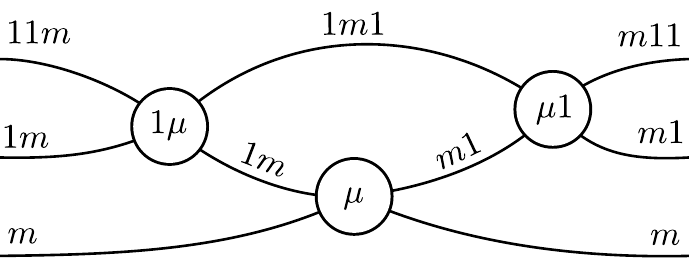}}

\begin{equation}\label{eqn:algebraassociativity}
\begin{tabular}{@{}ccc@{}}

\includegraphics[width=52.5mm]{Pictures/Preliminaries/Algebra/associativity1.pdf} & \raisebox{0.45\prelim}{$=$} &
\includegraphics[width=40mm]{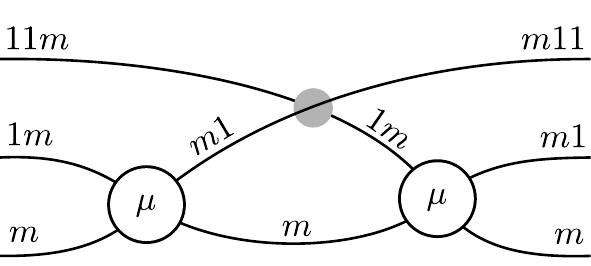},

\end{tabular}
\end{equation}

\begin{enumerate}
\item [b.] 
\end{enumerate}

\settoheight{\prelim}{\includegraphics[width=22.5mm]{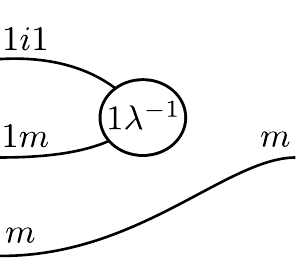}}

\begin{equation}\label{eqn:algebraunitality}
\begin{tabular}{@{}ccc@{}}

\includegraphics[width=22.5mm]{Pictures/Preliminaries/Algebra/unitality1.pdf} & \raisebox{0.45\prelim}{$=$} &

\includegraphics[width=37.5mm]{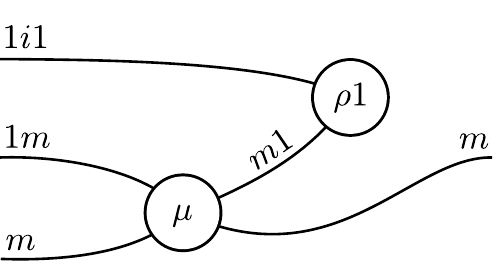}.

\end{tabular}
\end{equation}
\end{definition}

Let us now recall the definition of a right $A$-module in $\mathfrak{C}$ from \cite{D7}. Once, again the definition in a general monoidal 2-category may be found in \cite{D4}.

\begin{definition}\label{def:module}
A right $A$-module in $\mathfrak{C}$ consists of:
\begin{enumerate}
    \item An object $M$ of $\mathfrak{C}$;
    \item A 1-morphism $n^M:M\Box A\rightarrow M$;
    \item Two 2-isomorphisms
\end{enumerate}
\begin{center}
\begin{tabular}{@{}c c@{}}
$\begin{tikzcd}[sep=small]
MAA \arrow[dd, "1m"'] \arrow[rr, "n^M1"]    &  & MA \arrow[dd, "n^M"] \\
                                            &  &                      \\
MA \arrow[rr, "n^M"'] \arrow[rruu, Rightarrow, "\nu^M", shorten > = 2.5ex, shorten < = 2.5ex] &  & M,                   
\end{tikzcd}$

&

$\begin{tikzcd}[sep=small]
                                  &  & MA \arrow[rrdd, "n^M"] \arrow[dd, Rightarrow, "\rho^M", shorten > = 1ex, shorten < = 2ex] &  &   \\
                                  &  &                                             &  &   \\
M \arrow[rruu, "1i"] \arrow[rrrr,equal] &  & {}                                          &  & M,
\end{tikzcd}$
\end{tabular}
\end{center}

satisfying:

\begin{enumerate}
\item [a.] 
\end{enumerate}

\settoheight{\prelim}{\includegraphics[width=52.5mm]{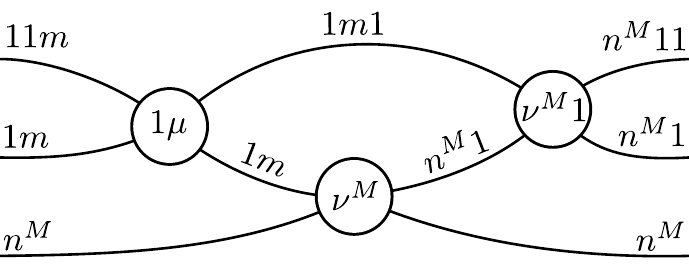}}

\begin{equation}\label{eqn:moduleassociativity}
\begin{tabular}{@{}ccc@{}}

\includegraphics[width=52.5mm]{Pictures/Preliminaries/Module/associativity1.pdf} & \raisebox{0.45\prelim}{$=$} &
\includegraphics[width=45mm]{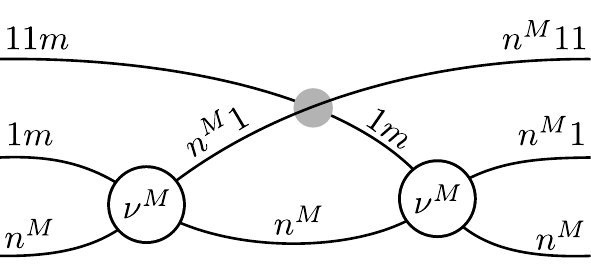},

\end{tabular}
\end{equation}

\begin{enumerate}
\item [b.]
\end{enumerate}

\settoheight{\prelim}{\includegraphics[width=22.5mm]{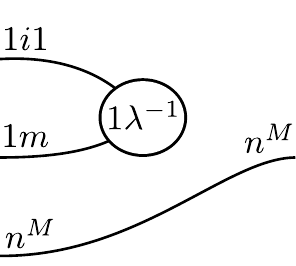}}

\begin{equation}\label{eqn:moduleunitality}
\begin{tabular}{@{}ccc@{}}

\includegraphics[width=22.5mm]{Pictures/Preliminaries/Module/unitality1.pdf} & \raisebox{0.45\prelim}{$=$} &

\includegraphics[width=37.5mm]{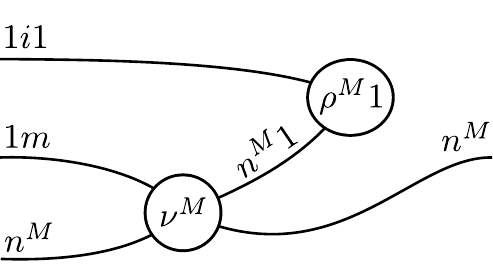}.

\end{tabular}
\end{equation}
\end{definition}

The definitions of a right $A$-module 1-morphism and that of a right $A$-module 2-morphism in $\mathfrak{C}$ may be found in \cite{D7}. These objects assemble into a 2-category as was proven in lemma 3.2.10 of \cite{D4}.

\begin{lemma}
Let $A$ be an algebra in a monoidal 2-category $\mathfrak{C}$. Right $A$-modules, right $A$-module 1-morphisms, and right $A$-module 2-morphisms form a 2-category, which we denote by $\mathbf{Mod}_{\mathfrak{C}}(A)$.
\end{lemma}

\subsection{Higher Condensations and Separable Algebras}

We now briefly review the notions of 2-condensations and 2-condensation monads. These notions were introduced in \cite{GJF} as the categorifications of the notions of split surjection and idempotent.

\begin{definition}
A 2-\define{condensation} in a $2$-category $\mathfrak{C}$ consists of two objects $A$ and $B$, together with two 1-morphisms $f : A \leftrightarrows B : g$, and two 2-morphisms $\phi:f\circ g\Rightarrow Id_B$ and $\gamma:Id_B\Rightarrow f\circ g$ , such that $\phi\cdot\gamma = Id_{Id_B}$.
\end{definition}
\noindent The data of 2-condensation as in the above definition induces a 2-condensation
monad on the object $A$.
\begin{definition}\label{def:2condmonad}
A 2-condensation monad in $\mathfrak{C}$ is an object $A$ together with a 1-morphism $e: A \to A$ and 2-morphisms $\mu: e\circ e \to e$ and $\delta: e \to e \circ e$, such that $\mu$ is associative, $\delta$ is coassociative, the Frobenius relations holds, and $\mu \cdot \delta = Id_{e}$.
\end{definition}
\noindent We say that a 2-condensation monad can be split, if it can be extended to a 2-condensation. There is also a categorification of the concept of idempotent complete 1-category. Before we review this definition, let us recall that a 2-category $\mathfrak{C}$ is locally idempotent complete if for all objects $A, B \in \mathfrak{C}$, the 1-category $\hom_\mathfrak{C}(A,B)$ is idempotent complete.
\begin{definition}
We say that a 2-category is \define{Karoubi complete} if it is locally idempotent complete, and every 2-condensation monad can be split.
\end{definition}
\noindent Physically, this means that any surface that arises as a condensation defect, i.e. a network of lower dimensional objects, is included in the 2-category.

The 2-category $\fC$ is locally finite semisimple if $\hom_\mathfrak{C}(A,B)$ is a finite semisimple $\bC$-linear 1-category (i.e.\ an abelian $\bC$-linear 1-category with finitely many isomorphism classes of simple object and in which every object decomposes as a finite direct sum of simple objects). We say that an object $A$ of $\fC$ is simple if the identity 1-morphism $Id_A$ is a simple object of the 1-category $\End_{\fC}(A)$.
\begin{definition}
A \define{finite semisimple} 2-category is a locally finite semisimple 2-category, that has adjoints for 1-morphisms, is Karoubi complete, has direct sums for objects, and has finitely many equivalence classes of objects.
\end{definition}
\noindent Finite semisimple 2-categories were introduced in \cite{DR}. We have recalled an equivalent version of their definition (see theorem 3.1.7 \cite{GJF}). Through proposition 1.4.5 of \cite{DR}, any object in a finite semisimple 2-category is the direct sum of finitely many simple objects, i.e. surfaces.

Let us recall the following definition from \cite{D2}. Thanks to section 2.2 of \cite{D2}, this is equivalent to the original definition given in \cite{DR}.
\begin{definition}
A multifusion 2-category is a finite semisimple rigid monoidal 2-category. A fusion 2-category is a multifusion 2-category whose monoidal unit is a simple object.
\end{definition}
\noindent Further, in a finite semisimple 2-category, two simple
objects that have a nonzero 1-morphism between
them are organized into the same component of $\mathfrak{C}$, denoted by $\pi_0(\mathfrak{C})$, due to the categorical Schur's lemma (see proposition 1.2.19 of \cite{DR}). In other words, $\pi_0(\mathfrak{C})$ only remembers objects up to condensation. We review the following definition from
 \cite{JFY}, due to its prevalence in section \ref{section:ExampleCats}:
\begin{definition}
A multifusion 2-category $\mathfrak{C}$ is \define{bosonic strongly fusion} if the braided fusion 1-category $\Omega \mathfrak{C} = \End_{\mathfrak{C}}(\mathbb{1}_\mathfrak{C}) $ is equivalent to $\Vec$. It is $\define{fermionic strongly fusion}$ if $\Omega \mathfrak{C} \simeq \SVec$.
\end{definition}
\noindent In such a 2-category $\fC$, the main result of \cite{JFY} shows that $\pi_0(\mathfrak{C})$ has grouplike fusion rules.

Definition \ref{def:2condmonad} has been categorified further in \cite{GJF} where the authors define an $n$-condensation monad for any $n$. Examples of 3-condensation monads are given by separable algebras in a monoidal 2-category as defined below. It is also convenient to introduce the notion of a rigid algebra, which can be traced back to \cite{G}. Rigid algebras are a weakening of separable algebras, and were first considered in the setting of fusion 2-categories in \cite{JFR}. We also point out that both of these definitions are thoroughly unpacked in section 2.1 of \cite{D7}.
\begin{definition}
An algebra $A$ in a monoidal 2-category $\mathfrak{C}$ is called rigid if the multiplication map $m:A\Box A\rightarrow A$ has a right adjoint $m^*$ as an $A$-$A$-bimodule 1-morphism. A rigid algebra $A$ in $\mathfrak{C}$ is called separable if the counit $\epsilon^m:m\circ m^*\Rightarrow Id_A$ witnessing that $m^*$ is right adjoint to $m$ as an $A$-$A$-bimodule 1-morphism has a section as an $A$-$A$-bimodule 2-morphism.
\end{definition}
We will see the separability property appear in the theorems in section \ref{section:BraidedSymmetric}. In fact, these results holds more generally for any 3-condensation monad. For later use, we also record the following result, which is given by combining together proposition 3.1.2 of \cite{D7} and corollary 2.2.3 of \cite{D5}.
\begin{proposition}\label{prop:module2categoryfinitesemisimple}
Let $A$ be a separable algebra in a fusion 2-category $\mathfrak{C}$. Then, the 2-category $\mathbf{Mod}_{\mathfrak{C}}(A)$ is a finite semisimple 2-category.
\end{proposition}

\begin{figure}[ht]
\centering
    \begin{tikzpicture}[thick]
    
        \def\Depth{5}
        \def\Height{3}
        \def\Width{3}
        \def\Sep{4}        
        
        \coordinate (O) at (0+1.5,0,0);
        \coordinate (A) at (1.5,\Width,0);
        \coordinate (B) at (1.5,\Width,\Height);
        \coordinate (C) at (0+1.5,0,\Height);
        \coordinate (D) at (\Depth,0,0);
        \coordinate (E) at (\Depth,\Width,0);
        \coordinate (F) at (\Depth,\Width,\Height);
        \coordinate (G) at (\Depth,0,\Height);
        \draw[black] (O) -- (D) -- (G) -- (C);
        \draw[black] (D) -- (E) -- (F) -- (G) -- cycle;
        \draw[black] (A) -- (E) -- (F) -- (B);
        \draw[midway] (\Depth/2,\Width-\Width/2,\Height/2) node {Category $\mathfrak{C}$};
        \draw[midway] (\Depth/2-2,\Width-\Width/2,\Height/2) node[scale=2] {\ldots};
        
        \draw[->,decorate,decoration={snake,amplitude=.4mm,segment length=2mm,post length=1mm}] (\Depth+1/2, \Width/2, \Height/2) -- (\Depth+\Sep-\Width/4,\Width/2,\Height/2) node[midway, above] {Condense $A$};
        
        \coordinate (OT) at (\Sep+\Depth+.25,0,0);
        \coordinate (AT) at (\Sep+\Depth+.25,\Width,0);
        \coordinate (BT) at (\Sep+\Depth,\Width,\Height);
        \coordinate (CT) at (\Sep+\Depth,0,\Height);
        \draw[black] (OT) -- (CT) -- (BT) -- (AT) -- cycle;
        \draw[midway] (\Sep+\Depth+.16, 1/2*\Width, 1/2*\Height) node[xslant = -0.5, yslant = 0.5]{$\Mod_{\mathfrak{C}}(A)$};
        \draw[midway] (\Sep+\Depth+2, 1/2*\Width, 1/2*\Height) node{$2\SVec$};
    \end{tikzpicture}
    \caption{This gives a three dimensional view of condensing the algebra $A$, taking place in a 2-category $\mathfrak{C}$. The resulting boundary is the category $\Mod_{\mathfrak{C}}(A)$, and $2\SVec$ represents the ``fermionic vacuum".  
    }
      \label{fig:condensebulk}
\end{figure}
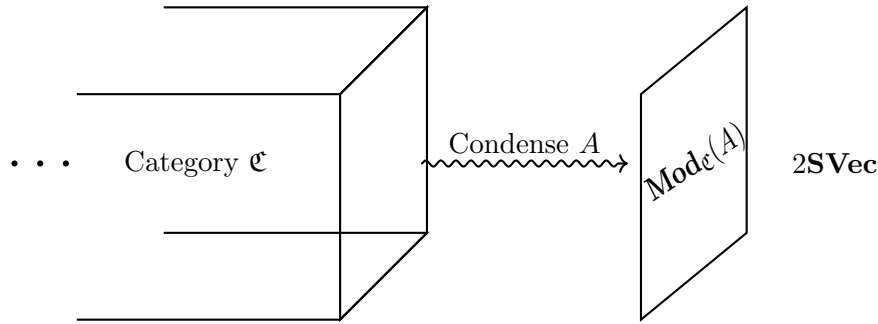

The physical picture for condensing surfaces in a 2-category involves finding some gapped boundary of the initial 2-category $\mathfrak{C}$, and then possibly triggering another condensation in order to map to $2\SVec$, see figure \ref{fig:condensebulk}. This bulk boundary point of view has been given the name of a ``quiche'', in \cite{FMT:2022}. The tensor unit of the boundary
can be identified with a separable algebra $A$ in $\mathfrak{C}$, and we denote it as $\Mod_\mathfrak{C}(A)$, the 2-category of $A$-modules in $\mathfrak{C}$. From this point of view, condensation along a specific direction of spacetime builds modules which usually causes the resulting 2-category to lose a level of monoidality, this is reflected in Theorems \ref{thm:first} and \ref{thm:second}. Theorem \ref{thm:third}, however, maintains the sylleptic property due to the extra condition of being in the symmetric center.
For a description of condensation in 1-categories where modules are explicitly built, see \cite{KO2014,Yu2021}. 

\subsection{Relative Tensor Product}\label{sub:tensor}
We now recall the definition of the relative tensor product over an algebra in a monoidal 2-category given in section 3 of \cite{D8}. These definitions will be important for the proofs of the main theorems in \S\ref{section:ModBraided} and \S\ref{section:ModSymmetric}. We also give sufficient criterion for the 2-category of bimodules over an algebra to carry a monoidal structure.

Let us now fix an algebra $A$ in a fusion 2-category $\mathfrak{C}$, together with $M$ a right $A$-module in $\mathfrak{C}$, and $N$ a left $A$-module in $\mathfrak{C}$ (for which we use the notations of appendix A of \cite{D7}).
We begin by defining $A$-balanced 1-morphisms and 2-morphisms out of the pair $(M,N)$.
\begin{definition}\label{def:balanced1morphism}
Let $C$ be an object of $\mathfrak{C}$. An $A$-balanced 1-morphism $(M,N)\rightarrow C$ consists of:
\begin{enumerate}
    \item A 1-morphism $f:M\Box N\rightarrow C$ in $\mathfrak{C}$;
    \item A 2-isomorphism
\end{enumerate}

$$\begin{tikzcd}[sep=small]
MAN \arrow[dd, "1l^N"'] \arrow[rr, "n^M1"]    &  & MN \arrow[dd, "f"] \\
                                            &  &                      \\
MN \arrow[rr, "f"'] \arrow[rruu, Rightarrow, "\tau^f", shorten > = 2.5ex, shorten < = 2.5ex] &  & C,                   
\end{tikzcd}$$

satisfying:

\begin{enumerate}
\item [a.]
\end{enumerate}

\newlength{\tensor}
\settoheight{\tensor}{\includegraphics[width=52.5mm]{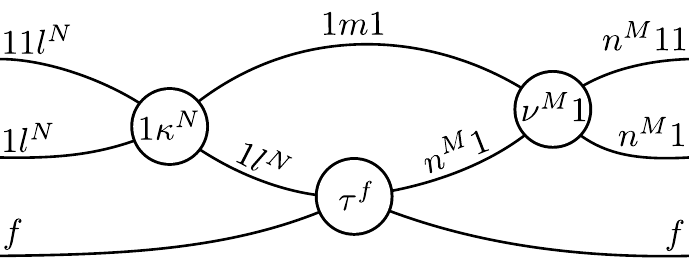}}

\begin{equation}\label{eqn:balancedassociativity}
\begin{tabular}{@{}ccc@{}}

\includegraphics[width=52.5mm]{Pictures/Preliminaries/Tensor/associativity1.pdf} & \raisebox{0.45\tensor}{$=$} &
\includegraphics[width=45mm]{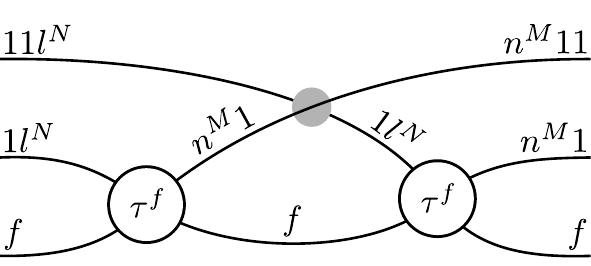},

\end{tabular}
\end{equation}

\begin{enumerate}
\item [b.]
\end{enumerate}

\settoheight{\tensor}{\includegraphics[width=22.5mm]{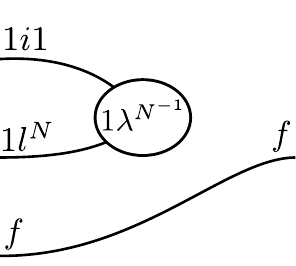}}

\begin{equation}\label{eqn:balancedunitality}
\begin{tabular}{@{}ccc@{}}

\includegraphics[width=22.5mm]{Pictures/Preliminaries/Tensor/unitality1.pdf} & \raisebox{0.45\tensor}{$=$} &

\includegraphics[width=37.5mm]{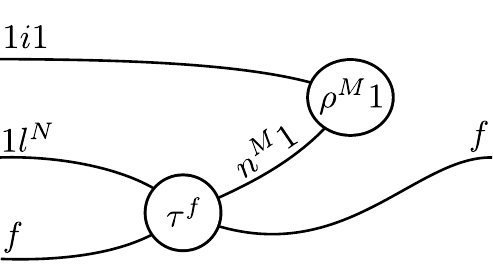}.

\end{tabular}
\end{equation}
\end{definition}

\begin{definition}\label{def:balanced2morphism}
Let $C$ be an object of $\mathfrak{C}$, and $f,g:(M,N)\rightarrow C$ be two $A$-balanced 1-morphisms. An $A$-balanced 2-morphism $f\Rightarrow g$ is a 2-morphism $\gamma:f\Rightarrow g$ in $\mathfrak{C}$ such that
\settoheight{\tensor}{\includegraphics[width=30mm]{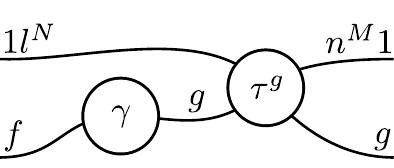}}

\begin{equation}\label{eqn:balance2morphism}
\begin{tabular}{@{}ccc@{}}

\includegraphics[width=30mm]{Pictures/Preliminaries/Tensor/2morphism1.pdf} & \raisebox{0.45\tensor}{$=$} &

\includegraphics[width=30mm]{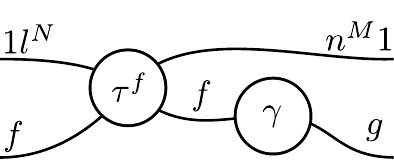}.

\end{tabular}
\end{equation}
\end{definition}

\begin{definition}\label{def:relativetensor}
The relative tensor product of $M$ and $N$ over $A$, if it exists, is an object $M\Box_A N$ of $\mathfrak{C}$ together with an $A$-balanced 1-morphism $t_A:(M,N)\rightarrow M\Box_A N$ satisfying the following 2-universal property:
\begin{enumerate}
    \item For every $A$-balanced 1-morphism $f:(M,N)\rightarrow C$, there exists a 1-morphism $\widetilde{f}:M\Box_A N\rightarrow C$ in $\mathfrak{C}$ and an $A$-balanced 2-isomorphism $\xi:\widetilde{f}\circ t_A\cong f$.
    \item For any 1-morphisms $g,h:M\Box_A N\rightarrow C$ in $\mathfrak{C}$, and any $A$-balanced 2-morphism $\gamma:g\circ t_A\Rightarrow h\circ t_A$, there exists a unique 2-morphism $\zeta:g\Rightarrow h$ such that $\zeta\circ t_A = \gamma$.
\end{enumerate}
\end{definition}

The following result was established in theorem 3.1.6 of \cite{D8}.

\begin{theorem}
Let $A$ be a separable algebra in a Karoubi complete monoidal 2-category $\mathfrak{C}$. Then, the relative tensor product of any right $A$-module $M$ and any left $A$-module $N$ exists.
\end{theorem}

Using this result, it was shown in theorem 3.2.8 of \cite{D8} that the relative tensor product over $A$ endows the 2-category $\mathbf{Bimod}_{\mathfrak{C}}(A)$ of $A$-$A$-bimodules in the Karoubi complete 2-category $\mathfrak{C}$ with a \textit{weak} monoidal structure. In particular, all the relevant structures were exhibited using the 2-universal property of the relative tensor product over multiple separable algebras.

\section{Braided and Symmetric Algebras}\label{section:BraidedSymmetric}

\subsection{Definitions}

Let $\mathfrak{B}$ be a semi-strict braided monoidal 2-category. The definition of a braided algebra in a braided monoidal 2-category, also called braided pseudo-monoid, can be traced back to \cite{DS}. Below we review this definition using the graphical calculus that we have previously introduced. We refer the reader to \cite{McC} for a version of this definition, resp. the next one, in a completely general braided, resp. sylleptic, monoidal 2-category.

\begin{definition}\label{def:braidedalgebra}
A braided algebra in $\mathfrak{B}$ consists of:
\begin{enumerate}
    \item An algebra $(B,m,i,\lambda, \mu, \rho)$ in $\mathfrak{B}$;
    \item A 2-isomorphisms
\end{enumerate}

$$\begin{tikzcd}[sep=small]
                                  &  & AA \arrow[rrdd, "m"] \arrow[dd, Rightarrow, "\beta", shorten > = 1ex, shorten < = 2ex] &  &   \\
                                  &  &                                             &  &   \\
AA \arrow[rruu, "b"] \arrow[rrrr, "m"'] &  & {}                                          &  & A,
\end{tikzcd}$$

satisfying:

\begin{enumerate}
\item [a.] 
\end{enumerate}

\newlength{\braid}

\settoheight{\braid}{\includegraphics[width=52.5mm]{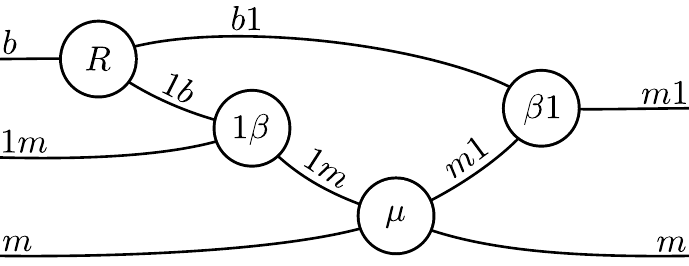}}

\begin{equation}\label{eqn:braidedalgebra1}
\begin{tabular}{@{}ccc@{}}

\includegraphics[width=52.5mm]{Pictures/BraidedAlgebrasProperties/BraidedAlgebra/braidedalgebra1.pdf} & \raisebox{0.45\braid}{$=$} &
\includegraphics[width=52.5mm]{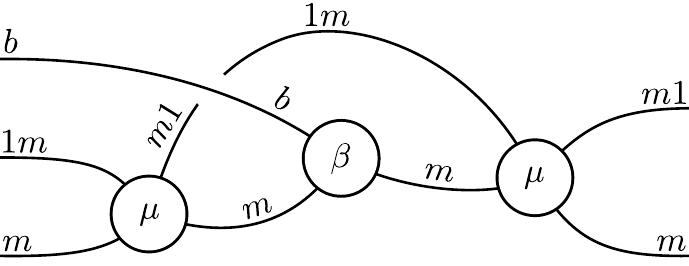},

\end{tabular}
\end{equation}

\begin{enumerate}
\item [b.] 
\end{enumerate}

\settoheight{\braid}{\includegraphics[width=52.5mm]{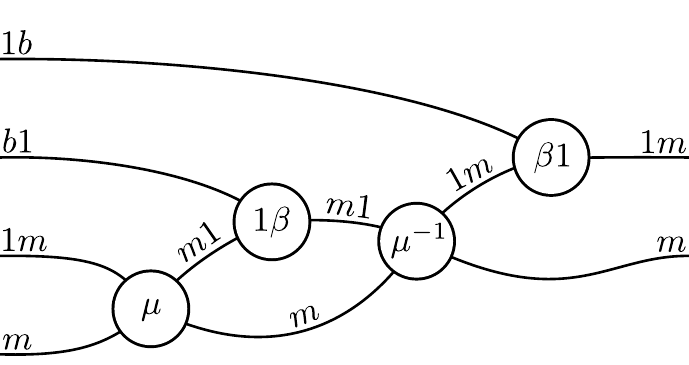}}

\begin{equation}\label{eqn:braidedalgebra2}
\begin{tabular}{@{}ccc@{}}

\includegraphics[width=52.5mm]{Pictures/BraidedAlgebrasProperties/BraidedAlgebra/braidedalgebra3.pdf} & \raisebox{0.45\braid}{$=$} &

\includegraphics[width=52.5mm]{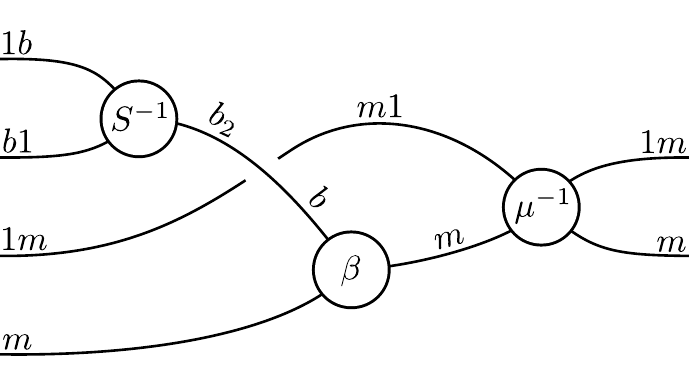}.

\end{tabular}
\end{equation}

\begin{enumerate}
\item [c.] 
\end{enumerate}

\settoheight{\braid}{\includegraphics[width=30mm]{Pictures/BraidedAlgebrasProperties/BraidedAlgebra/braidedalgebra4.pdf}}

\begin{equation}\label{eqn:braidedalgebra3}
\begin{tabular}{@{}ccc@{}}

\includegraphics[width=30mm]{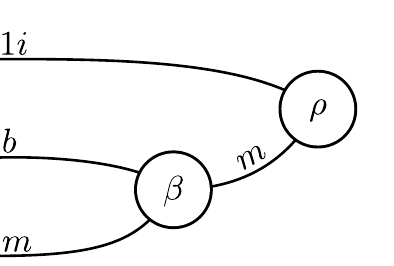} & \raisebox{0.45\braid}{$=$} &

\includegraphics[width=30mm]{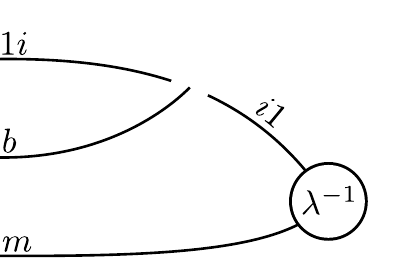}.

\end{tabular}
\end{equation}
\end{definition}

Let $\mathfrak{S}$ be a semi-strict sylleptic monoidal 2-category. The definition of a symmetric algebra in $\mathfrak{S}$, also called symmetric pseudo-monoid, first appeared in \cite{DS}. We review this definition using our graphical calculus.

\begin{definition}
A symmetric algebra in $\mathfrak{S}$ is a braided algebra $(B,m,i,\lambda, \mu, \rho, \beta)$ such that 

\settoheight{\braid}{\includegraphics[width=45mm]{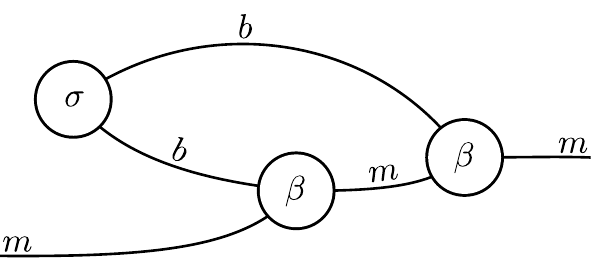}}

\begin{equation}\label{eqn:symmetricalgebra}
\begin{tabular}{@{}cc@{}}

\includegraphics[width=45mm]{Pictures/BraidedAlgebrasProperties/BraidedAlgebra/symmetricalgebra.pdf}

&\raisebox{0.45\braid}{$= Id_m.$} 
\end{tabular}
\end{equation}

\end{definition}

\begin{example}
Braided algebras in the symmetric fusion 2-category $2\Vec$ are exactly braided monoidal finite semisimple 1-categories. Symmetric algebras in the symmetric fusion 2-category $2\Vec$ are exactly symmetric monoidal finite semisimple 1-categories.
\end{example}

\subsection{The 2-Category of Modules over a Braided Algebra.}\label{section:ModBraided}

As before, we take $\mathfrak{B}$ to be a semi-strict braided monoidal 2-category. Furthermore, we will assume throughout that $\mathfrak{B}$ is a Karoubi complete 2-category.

\begin{lemma}\label{lem:modulesub2category}
Let $B$ a braided algebra in $\mathfrak{B}$. There is a 2-functor $$Ind^+:\mathbf{Mod}_{\mathfrak{B}}(B)\rightarrow \mathbf{Bimod}_{\mathfrak{B}}(B),$$ which is fully faithful on 2-morphisms.
\end{lemma}
\begin{proof}
Let $M$ be a right $B$-module. The underlying right $B$-module of $Ind^+(B)$ is given by $B$. In the notations of \cite{D7}, the left $B$-module structure on $Ind^+(B)$ is given by the 1-morphism $$l^M:B\Box M\xrightarrow{b}M\Box B\xrightarrow{n^M}M$$ together with the 2-isomorphisms $$\raisebox{15pt}{\settoheight{\braid}{\includegraphics[width=30mm]{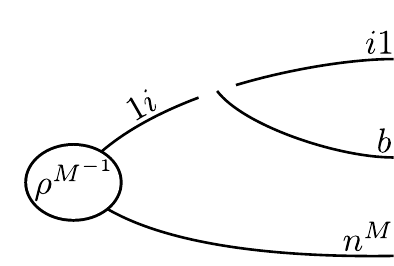}}
\raisebox{0.45\braid}{$\lambda^M:=\ $}
\includegraphics[width=30mm]{Pictures/BraidedAlgebrasProperties/ModuleBraidedAlgebra/lambdaM.pdf},}\ \ \ \settoheight{\braid}{\includegraphics[width=52.5mm]{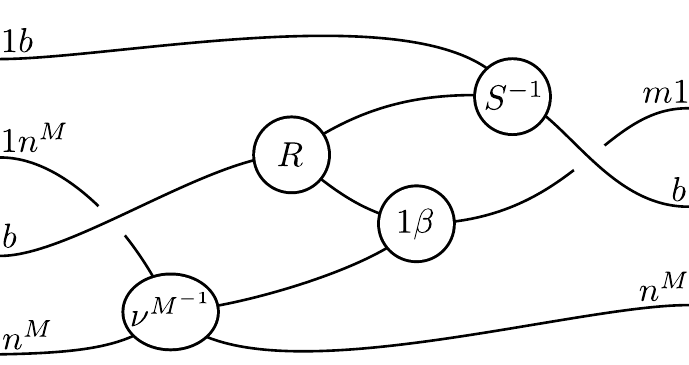}}
\raisebox{0.45\braid}{$\kappa^M:=\ $}
\includegraphics[width=52.5mm]{Pictures/BraidedAlgebrasProperties/ModuleBraidedAlgebra/kappaM.pdf}.$$ Further, the compatibility between the left and the right actions is given by the 2-isomorphism $$\settoheight{\braid}{\includegraphics[width=52.5mm]{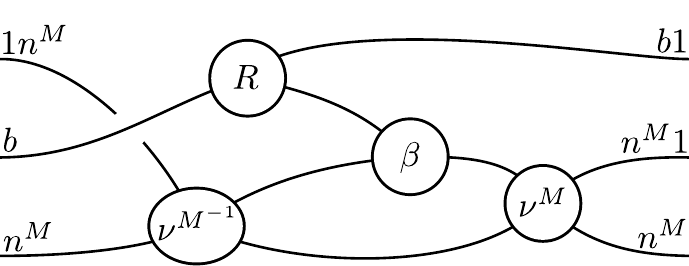}}
\raisebox{0.45\braid}{$\beta^M:=\ $}
\includegraphics[width=52.5mm]{Pictures/BraidedAlgebrasProperties/ModuleBraidedAlgebra/betaM.pdf}.$$

Given a right $B$-module 1-morphism $f:M\rightarrow N$, the underlying right $B$-module 1-morphism of the $B$-$B$-module 1-morphism $Ind^+(f)$ is $f$. Its left $B$-module structure is given by $$\settoheight{\braid}{\includegraphics[width=30mm]{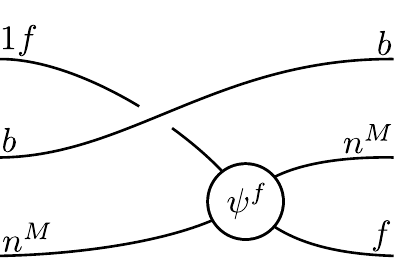}}
\raisebox{0.45\braid}{$\chi^f:=\ $}
\includegraphics[width=30mm]{Pictures/BraidedAlgebrasProperties/ModuleBraidedAlgebra/chif.pdf}.$$ Given a right $B$-module 2-morphism $\gamma:f\Rightarrow g$, it is easy to check that $\gamma$ is a $B$-$B$-bimodule 2-morphism $Ind^+(f)\Rightarrow Ind^+(g)$, so that we can set $Ind^+(\gamma)=\gamma$. It follows readily from the definitions that $Ind^+$ defines a strict 2-functor. Moreover, note that $Ind^+$ is fully faithful on 2-morphisms by construction.\end{proof}

\begin{rem}
When constructing the 2-functor $Ind^+$, we have used the braiding $b$ of $\mathfrak{B}$. Instead, we could have used its adjoint equivalence $b^{\bullet}$, and so doing obtained a 2-functor $Ind^-:\mathbf{Mod}_{\mathfrak{B}}(B)\rightarrow \mathbf{Bimod}_{\mathfrak{B}}(B)$.
\end{rem}

\begin{proposition}\label{prop:braidedalgebramodules}
Let $B$ a braided separable algebra in $\mathfrak{B}$. Then, $\mathbf{Mod}_{\mathfrak{B}}(B)$ is a monoidal 2-category with monoidal unit $B$.
\end{proposition}
\begin{proof}
Thanks to lemma \ref{lem:modulesub2category}, we can view $\mathbf{Mod}_{\mathfrak{B}}(B)$ as a sub-2-category of $\mathbf{Bimod}_{\mathfrak{B}}(B)$. For convenience, we will assume that this sub-2-category is replete. Now, as was recalled in section \ref{sub:tensor} above, the monoidal structure of $\mathbf{Bimod}_{\mathfrak{B}}(B)$ is given by the {relative tensor product} $\Box_B$, which is defined using the 2-universal property reviewed in definition \ref{def:relativetensor}.

Given $M$ and $N$ two right $B$-modules, we want to show that the $B$-$B$-bimodule $M\Box_B N$ is actually an object of the sub-2-category $\mathbf{Mod}_{\mathfrak{B}}(B)$. In order to prove this, we need to unfold the definition of the left $B$-module structure on $M\Box_B N$. Let us write $t:M\Box N\rightarrow M\Box_B N$, together with $\tau^t:t\circ (M\Box l^N)\cong t\circ (n^M\Box N)$, for the 2-universal $B$-balanced 1-morphism as in definition \ref{def:relativetensor}. Furthermore, note that for any $C$ in $\mathfrak{B}$, $C\Box t$ equipped with $C\Box \tau^t$ is a 2-universal $B$-balanced 1-morphism. By remark 3.2.3 of \cite{D8}, the 1-morphism $l^{M\Box_BN}:B\Box (M\Box_B N)\rightarrow M\Box_B N$ is induced by the 2-universal property of $B\Box t$ applied to the solid arrow diagram $$\begin{tikzcd}[sep=small]
B\Box M\Box N \arrow[dd, "(n^M\circ b)1"'] \arrow[rr, "1t"] &  & B\Box (M\Box_B N) \arrow[dd, "l^{M\Box_BN}", dotted]  \\
                                                            &  &                                                                     \\
M\Box N \arrow[rr, "t"']  \arrow[rruu, Rightarrow, "\upsilon^l", shorten > = 3ex, shorten < = 3ex]                                  &  & M\Box_B N ,                                                         
\end{tikzcd}$$ where the left bottom composite 1-morphism is equipped with the obvious $B$-balanced structure. The 1-morphism $n^{N\Box_B M}:(M\Box_B N)\Box B\rightarrow M\Box_B N$ is defined similarly. But, the 2-isomorphism $$(t\circ 1n^N \circ R^{-1})\cdot(\tau^{t^{-1}}\circ b1):t\circ (n^M\Box N) \circ (b\Box N)\cong t\circ (M\Box n^N)\circ b$$ is $B$-balanced. Thanks to the 2-universal property of the relative tensor product, this means that there exists a 2-isomorphism $\theta:l^{M\Box_BN}\cong n^{M\Box_BN}\circ b$. Furthermore, it also follows from the 2-universal property that $\theta$ promotes the identity right $B$-module 1-morphism on $M\Box_BN$ to a $B$-$B$-bimodule 1-equivalence from $M\Box_BN$ to $Ind^+(M\Box_BN)$. This proves that the objects of $\mathbf{Mod}_{\mathfrak{B}}(B)$ are closed under $\Box_B$. A similar argument shows that the 1-morphisms of $\mathbf{Mod}_{\mathfrak{B}}(B)$ are closed under $\Box_B$, which concludes the proof.
\end{proof}

\begin{rem}\label{rem:generaltensor}
We emphasize that $\mathbf{Mod}_{\mathfrak{B}}(B)$ is not a braided 2-category in general, as can be seen from example \ref{ex:braidedfusion} below. Further, we also note that our proof of proposition \ref{prop:braidedalgebramodules} only used the existence of the relative tensor product over $B$ for any $B$-$B$-bimodules in $\mathfrak{B}$. We refer the reader to remark 3.2.11 of \cite{D8} for a more thorough discussion. An analogous comment can be made with regards to lemma and \ref{lem:modulesub2category} above and lemma \ref{lem:freemodulemonoidal} below.
\end{rem}

In order to prove our next theorem, we need the following technical lemma.

\begin{lemma}\label{lem:freemodulemonoidal}
The 2-functor $F:\mathfrak{B}\rightarrow \mathbf{Mod}_{\mathfrak{B}}(B)$ given by sending the object $C$ in $\mathfrak{B}$ to $C\Box B$ with its canonical right $B$-module structure is monoidal.
\end{lemma}
\begin{proof}
Let $C$ and $D$ be two objects of $\mathfrak{B}$. Firstly, note that $C\Box D \Box B$ satisfies the 2-universal property of $(C\Box B)\Box_B (D\Box B)$ in $\mathbf{Bimod}_{\mathfrak{B}}(B)$. More precisely, the $B$-$B$-bimodule 1-morphism $$u_{C,D}:C\Box B\Box D\Box B\xrightarrow{1b1} C\Box D\Box B\Box B\xrightarrow{11m}C\Box D\Box B$$ admits a canonical $B$-balanced structure given by
$$\settoheight{\braid}{\includegraphics[width=52.5mm]{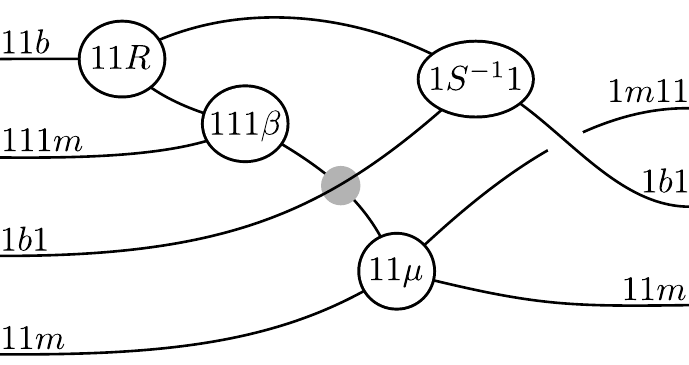}}
\raisebox{0.45\braid}{$\tau^{u_{C,D}}:=\ $}
\includegraphics[width=52.5mm]{Pictures/BraidedAlgebrasProperties/ModuleBraidedAlgebra/tauuCD.pdf}$$ and satisfies the conditions of definition \ref{def:relativetensor}. In particular, this yields $B$-$B$-bimodule 1-equivalences $e_{C,D}:(C\Box B)\Box_B (D\Box B)\simeq (C\Box D)\Box B$ for every $C$ and $D$ in $\mathfrak{B}$ together with a $B$-balanced $B$-$B$-bimodule 2-ismorphism $\zeta_{C,D}$ as in the following diagram: $$\begin{tikzcd}
C\Box B\Box D\Box B \arrow[rr, "u_{C,D}"] \arrow[rd, "t_{CB,DB}"'] & {} \arrow[d,Rightarrow, "\zeta_{C,D}"]             & C\Box D\Box B \\
                                      & (C\Box B)\Box_B (D\Box B).  \arrow[ru, "e_{C,D}"'] &    
\end{tikzcd}$$ Secondly, observe that for any two 1-morphisms $f:C\rightarrow E$ and $g:D\rightarrow F$ in $\mathfrak{B}$, the $B$-$B$-bimodule 2-isomorphism $$\settoheight{\braid}{\includegraphics[width=45mm]{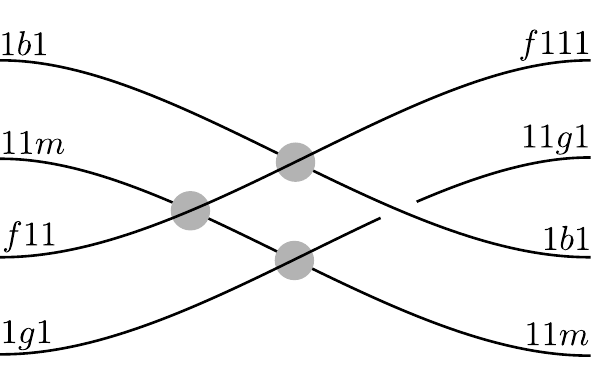}}
\raisebox{0.45\braid}{$\upsilon_{f,g}:=\ $}
\includegraphics[width=45mm]{Pictures/BraidedAlgebrasProperties/ModuleBraidedAlgebra/upsilonfg.pdf}$$ is $B$-balanced. Thus, thanks to the 2-universal property of the relative tensor product, we can use the 2-isomorphisms $\upsilon_{f,g}$ to promote the collection of the $B$-$B$-bimodule 1-equivalences $e_{C,D}$ for varying $C$ and $D$ to a 2-natural equivalence $e$.

Using the 2-universal property of the relative tensor product repeatedly (together with the variants over multiple algebras considered in section 3.2 of \cite{D8}), one constructs the remaining data necessary to endow $F$ with a monoidal structure, and prove that they satisfy the relevant axioms from definition 2.5 of \cite{SP}.
\end{proof}

\begin{theorem}\label{thm:multifusion2cat}
Let $\mathfrak{B}$ be a braided multifusion 2-category, and $B$ a braided separable algebra in $\mathfrak{B}$. Then, $\mathbf{Mod}_{\mathfrak{B}}(B)$ is a multifusion 2-category.
\end{theorem}

\begin{proof}
The 2-category $\mathbf{Mod}_{\mathfrak{B}}(B)$ is finite semisimple thanks to proposition \ref{prop:module2categoryfinitesemisimple}. Further, we have shown in proposition \ref{prop:braidedalgebramodules} that it admits a monoidal structure. It therefore only remains to prove that it has duals. But, as $\mathfrak{B}$ is a multifusion 2-category, it has right and left duals. In particular, every object in the image of $F:\mathfrak{B}\rightarrow \mathbf{Mod}_{\mathfrak{B}}(B)$ has a right and a left dual. But, it was shown in lemma 3.1.1 of \cite{D7} that every right $B$-module $M$ is the splitting of a 2-condensation monad (in $\mathbf{Mod}_{\mathfrak{B}}(B)$) supported on $M\Box B = F(M)$. Thence, it follows from lemma 5.5 of \cite{D3} that $M$ has a right and a left dual, and thereby concludes the proof.
\end{proof}

Following section 5.2 of \cite{D8}, we say that a separable algebra $B$ is connected if its unit 1-morphism $i:I\rightarrow B$ is simple. Under the equivalence $$Hom_{B}(B,B)\simeq Hom_{\mathfrak{B}}(I,B)$$ of lemma 3.2.13 of \cite{D4}, we have $Id_B\mapsto i$. Thus, $B$ is a simple right $B$-module if and only if $B$ is a connected algebra. Combined with the above theorem, this yields the following corollary.

\begin{corollary}
Let $\mathfrak{B}$ be a braided multifusion 2-category, and $B$ a connected braided separable algebra in $\mathfrak{B}$. Then, $\mathbf{Mod}_{\mathfrak{B}}(B)$ is a fusion 2-category.
\end{corollary}

\begin{example}\label{ex:braidedfusion}
Let $\mathcal{B}$ be a braided multifusion 1-categories, that is a braided separable algebra in $2\Vec$. Then, $\Mod_{2\Vec}(\mathcal{B})=\Mod(\mathcal{B})$ is the multifusion 2-category of finite semisimple right $\mathcal{B}$-module 1-categories with monoidal structure given by $\boxtimes_{\mathcal{B}}$ the relative Deligne tensor over $\mathcal{B}$. The braided separable algebra $\mathcal{B}$ is braided if and only if $\mathcal{B}$ is a fusion 1-category, in which case $\Mod(\mathcal{B})$ is a fusion 2-category. Finally, we note that it follows from a slight variant of proposition 2.4.7 of \cite{D2} that $\Mod(\mathcal{B})$ is braided if and only if $\mathcal{B}$ is symmetric.
\end{example}

\subsection{The 2-Category of Modules over a Symmetric Algebra}\label{section:ModSymmetric}
In this section we give sufficient conditions for the 2-category of modules over a braided algebra to be itself braided. We also explain when the 2-category of modules is sylleptic or symmetric.

\begin{theorem}\label{thm:SymmetricAlgebraModule}
Let $\mathfrak{S}$ be a Karoubi complete sylleptic monoidal 2-category, and $B$ a symmetric separable algebra in $\mathfrak{S}$. Then, $\mathbf{Mod}_{\mathfrak{S}}(B)$ is a braided monoidal 2-category.
\end{theorem}
\begin{proof}
Without loss of generality, we may assume that $\mathfrak{S}$ is semi-strict. Our first task is to endow the monoidal 2-category $\mathbf{Mod}_{\mathfrak{S}}(B)$ with a braiding $\widetilde{b}$. To this end, let $M$ and $N$ be two right $B$-modules, and write $$t_{M,N}:M\Box N\rightarrow M\Box_B N\ \textrm{and}\ t_{N,M}:N\Box M\rightarrow N\Box_B M$$ for the 2-universal $B$-balanced right $B$-module 1-morphisms with structure 2-isomorphisms $\tau^t$. We claim that the 1-morphism $t_{N,M}\circ b_{M,N}:M\Box N\rightarrow N\Box_B M$ in $\mathfrak{S}$ can be upgraded to a $B$-balanced right $B$-module 1-morphism. Namely, the $B$-balanced structure is given by the 2-isomorphism $$\settoheight{\braid}{\includegraphics[width=90mm]{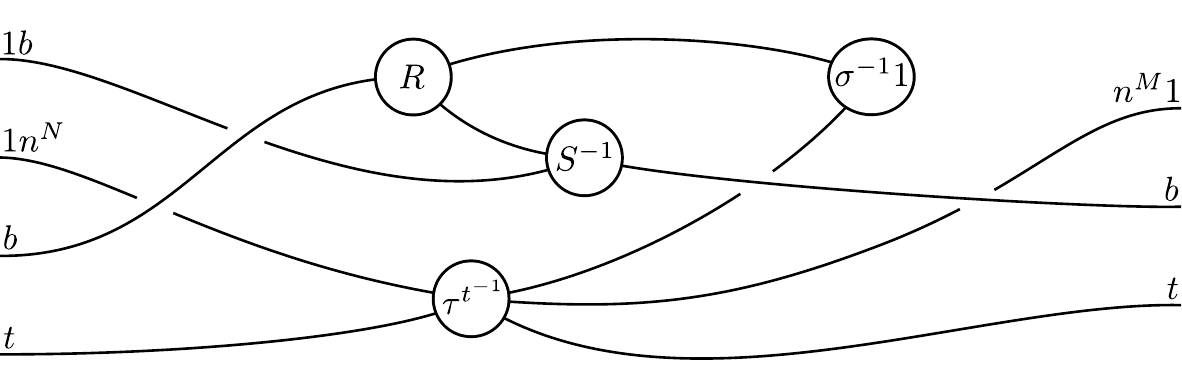}}
\raisebox{0.45\braid}{$\tau^{t\circ b}:=\ $}
\includegraphics[width=90mm]{Pictures/BraidedAlgebrasProperties/ModulesSymmetricAlgebras/tautb.pdf}.$$ In order to check that $\tau^{t\circ b}$ satisfies axiom a of definition \ref{def:relativetensor}, we use the diagrams depicted in appendix \ref{subsub:proofSymmetricAlgebraModule}. Figure \ref{fig:taubbalanced1} depicts the right hand-side of equation (\ref{eqn:balancedassociativity}). By moving the indicated coupons to the top along the blue arrows, we arrive at figure \ref{fig:taubbalanced2}. Then, using equation (\ref{eqn:balancedassociativity}) for $\tau^t$ on the blue coupons brings us to figure \ref{fig:taubbalanced3}. At this point, we use the definition of $\kappa^M$ given in the proof of lemma \ref{lem:modulesub2category} on the blue coupon, which leads us to contemplate figure \ref{fig:taubbalanced4}. Moving the coupon labeled $11\beta^{-1}$ to the left along the blue arrow, and then using equation (\ref{eqn:braidingaxiom4}) on the green coupons brings us to figure \ref{fig:taubbalanced5}. We arrive at figure \ref{fig:taubbalanced6} by moving the blue coupons to the left along the blue arrows. Moving the coupon labeled $1R^{-1}$ to the right along the blue arrow and then applying equation (\ref{eqn:syllepsisaxiom2}) on the green coupons bring us to figure \ref{fig:taubbalanced7}. By moving the coupon labeled $1S$ to the right along the blue arrow and then make use of equation (\ref{eqn:syllepsisaxiom1}) on the green coupons bring us to figure \ref{fig:taubbalanced8}. Using equation (\ref{eqn:braidingaxiom3}) on the blue coupons, we arrive at figure \ref{fig:taubbalanced9}. We obtain figure \ref{fig:taubbalanced10} by applying equation (\ref{eqn:braidingaxiom1}) on the blue coupons, using equation (\ref{eqn:braidingaxiom2}) on the green coupons, and moving the coupon labeled $1\beta^{-1}1$ to the right along the red arrow. Then, using equation (\ref{eqn:braidingaxiom4}) on the blue coupons and equation (\ref{eqn:symmetricalgebra}) on the green coupons, we arrive at figure \ref{fig:taubbalanced11}. Finally, we get to figure \ref{fig:taubbalanced12}, which depicts the left hand-side of equation (\ref{eqn:balancedassociativity}), by moving the coupon labeled $R$ to the right along the blue arrow and the coupon labeled $\beta^{-1}11$ to the left along the green arrow. Furthermore, equation (\ref{eqn:balancedunitality}) for $\tau^{t\circ b}$ follows from equation (\ref{eqn:balancedunitality}) for $\tau^t$ together with the fact that $R$, $S$, $\sigma$ are modifications, combined with axiom f of definition \ref{subsub:braidedmonoidal2cat} and axiom c of definition \ref{subsub:syllepticmonoidal2cat}.

Moreover, the right $B$-module structure on $t_{N,M}\circ b_{M,N}$ is given by the 2-isomorphism $$\settoheight{\braid}{\includegraphics[width=52.5mm]{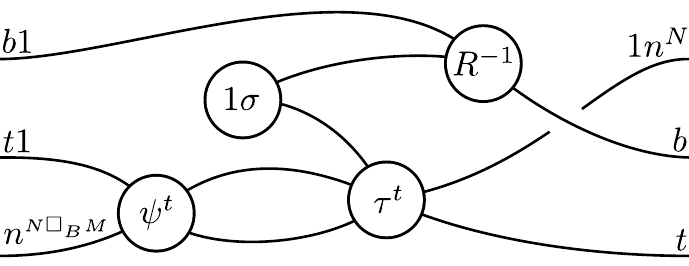}}
\raisebox{0.45\braid}{$\psi^{t\circ b}:=\ $}
\includegraphics[width=52.5mm]{Pictures/BraidedAlgebrasProperties/ModulesSymmetricAlgebras/psitb.pdf}.$$ Thus, by the 2-universal property of $t_{M,N}$, the solid arrow diagram below can be filled by a $B$-balanced right $B$-module 2-isomorphism $\xi_{M,N}$: $$\begin{tikzcd}[sep=small]
M\Box N \arrow[dd, "b_{M,N}"'] \arrow[rr, "t_{M,N}"] &  & M\Box_B N \arrow[dd, "\widetilde{b}_{M,N}", dotted]  \\
                                                            &  &                                                                     \\
N\Box M \arrow[rr, "t_{N,M}"'] \arrow[rruu, Rightarrow, "\xi_{M,N}", shorten > = 3ex, shorten < = 3ex]                                   &  & N\Box_B M.                                                         
\end{tikzcd}$$ Furthermore, as $b_{M,N}$ is a 1-equivalence, the 2-universal property implies that the 1-morphism $\widetilde{b}_{M,N}$ is also an equivalence. Using the 2-universal property of the relative tensor product over $B$ again, we find that the collection of the 1-equivalences $\widetilde{b}_{M,N}$ assembles to form a 2-natural equivalence $\widetilde{b}$. We upgrade $\widetilde{b}$ to an adjoint 2-natural equivalence by appealing to the 2-universal property.

We also have to construct invertible modifications $\widetilde{R}$ and $\widetilde{S}$ witnessing the coherence of the braiding $\widetilde{b}$ on $\mathbf{Mod}_{\mathfrak{S}}(B)$. As the monoidal structure on $\mathbf{Mod}_{\mathfrak{S}}(B)$ is not strict cubical, we need to use the fully weak definition of these modifications given in figure 2.3 of \cite{SP}. Let $M$, $N$, and $P$ be three right $B$-modules, in order to construct the right $B$-module 2-isomorphism $\widetilde{R}_{M,N,P}$ we use the 2-universal property of the relative tensor product over two algebras following definition 3.2.6 of \cite{D8}. More precisely, let us consider the 3-dimensional commutative diagram whose back and front are depicted below: $$\begin{tikzcd}[sep=small]
                                                                        & N\Box M\Box P \arrow[rr,equal] \arrow[dd]             &                                & N\Box M\Box P \arrow[dd] \arrow[rd, "1b"]         &                          \\
M\Box N\Box P \arrow[ru, "b1"] \arrow[dd]                               &                                                 &                                &                                                   & N\Box P\Box M \arrow[dd] \\
                                                                        & (N\Box_B M)\Box_B P \arrow[rr, "\alpha"']       & {} \arrow[dd, "\exists !\widetilde{R}", Rightarrow, dashed, shorten > = 3ex, shorten < = 3ex] & N\Box_B (M\Box_B P) \arrow[rd, "1\widetilde{b}"'] &                          \\
(M\Box_B N)\Box_B P \arrow[ru, "\widetilde{b}1"'] \arrow[rd, "\alpha"'] &                                                 &                                &                                                   & {N\Box_B (P\Box_B M),}   \\
                                                                        & M\Box_B(N\Box_B P) \arrow[rr, "\widetilde{b}"'] & {}                             & (N\Box_B P)\Box_B M \arrow[ru, "\alpha"']         &                         
\end{tikzcd}$$

$$\begin{tikzcd}[sep=small]
                                                     & N\Box M\Box P \arrow[rr,equal]                         & {} \arrow[dd, Rightarrow, "R", shorten > = 3ex, shorten < = 3ex] & N\Box M\Box P \arrow[rd, "1b"]            &                          \\
M\Box N\Box P \arrow[ru, "b1"] \arrow[rd,equal] \arrow[dd] &                                                  &                    &                                           & N\Box P\Box M \arrow[dd] \\
                                                     & M\Box N\Box P \arrow[rr, "b"'] \arrow[dd]        & {}                 & N\Box P\Box M \arrow[ru,equal] \arrow[dd]       &                          \\
(M\Box_B N)\Box_B P \arrow[rd, "\alpha"']            &                                                  &                    &                                           & N\Box_B (P\Box_B M).      \\
                                                     & M\Box_B (N\Box_B P) \arrow[rr, "\widetilde{b}"'] &                    & (N\Box_B P)\Box_B M \arrow[ru, "\alpha"'] &                         
\end{tikzcd}
$$

\noindent All the vertical 1-morphisms are 2-universal $(B,B)$-balanced right $B$-module 1-morphisms, and all the square faces are filled by $(B,B)$-balanced right $B$-module 2-isomorphisms thanks to either the proof of lemma 3.2.7 of \cite{D8} or the construction of $\widetilde{b}$ given above. Thus, thanks to the 2-universal property of the relative tensor product, there exists a unique right $B$-module 2-isomorphism $\widetilde{R}$ such that the whole 3-dimensional prism is commutative. Furthermore, the collection of these assignments assemble into an invertible modification as can been seen using the 2-universal property of the relative tensor product over two algebras. The invertible modification $\widetilde{S}$ is constructed similarly.

Finally, one has to check that $\widetilde{R}$ and $\widetilde{S}$ together with the modifications supplied by the monoidal structure of $\mathbf{Mod}_{\mathfrak{S}}(B)$ satisfy the equations given in figures C.7 through C.14 of \cite{SP} hold. This follows from the 2-universal property of the relative tensor product over three and four algebras explained in the proof of theorem 3.2.8 of \cite{D8}.
\end{proof}

\begin{proposition}\label{prop:freemodulebraided}
Let $\mathfrak{S}$ be a Karoubi complete sylleptic monoidal 2-category, and $B$ a symmetric separable algebra in $\mathfrak{S}$. Then, the monoidal 2-functor $F:\mathfrak{S}\rightarrow \mathbf{Mod}_{\mathfrak{S}}(B)$ of lemma \ref{lem:freemodulemonoidal} is braided.
\end{proposition}
\begin{proof}
Let $C$ and $D$ be any objects of $\mathfrak{S}$. Using the notations of lemma \ref{lem:freemodulemonoidal} and theorem \ref{thm:SymmetricAlgebraModule}, we can consider the following diagram: \begin{equation}\label{eqn:definitionepsilon}\begin{tikzcd}[sep=small]
                                                                         &  & {} \arrow[d,Rightarrow, "\:\zeta"]                                       &  &                      \\
CBDB \arrow[dd, "b_2"'] \arrow[rr, "t"] \arrow[rrrr, "u", bend left]     &  & (CB)\Box_B(DB) \arrow[dd, "\widetilde{b}"] \arrow[rr, "e"] &  & CDB \arrow[dd, "b1"] \\
                                                                         &  &                                                             &  &                      \\
DBCB \arrow[rr, "t"'] \arrow[rrrr, "u"', bend right] \arrow[rruu,Rightarrow, "\xi", shorten > = 4ex, shorten < = 4ex] &  & (DB)\Box_B(CB) \arrow[rr, "e"'] \arrow[d,Rightarrow, "\:\zeta^{-1}"]\arrow[rruu,Rightarrow,dashed, "\exists ! \epsilon", shorten > = 4ex, shorten < = 4ex]           &  & DCB.                  \\
                                                                         &  & {}                                                          &  &                     
\end{tikzcd}\end{equation} Further, the outer square can be filled using the $B$-balanced right $B$-module 2-isomorphism $\varsigma$ given by: $$\settoheight{\braid}{\includegraphics[width=60mm]{Pictures/BraidedAlgebrasProperties/ModuleBraidedAlgebra/tauuCD.pdf}}
\raisebox{0.45\braid}{$\varsigma:=\ $}
\includegraphics[width=60mm]{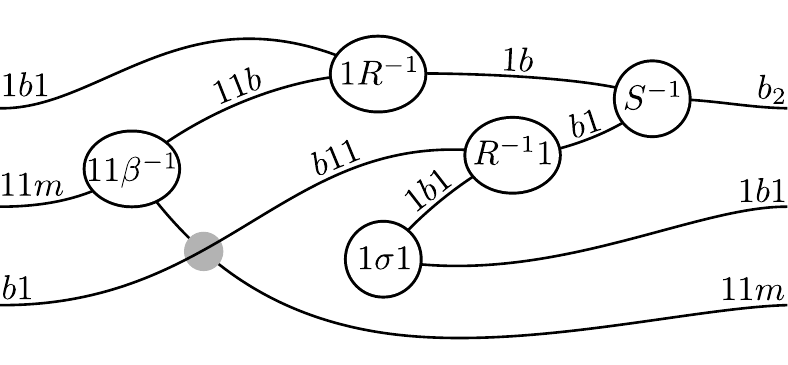}.$$ Thus, thanks to the 2-universal property of the relative tensor product over $B$, the right hand-side square of the commutative diagram (\ref{eqn:definitionepsilon}) can be filled by a right $B$-module 2-isomorphism $\epsilon$ such that its full composite is equal to $\varsigma$. Further, it follows from the same 2-universal property that the collection of these 2-isomorphism defines an invertible modification. Finally, one checks that the axioms of definition 2.5 of \cite{SP} hold for $\epsilon$ using the 2-universal property of the relative tensor product over one and two algebras.
\end{proof}

Note that $\mathbf{Mod}_{\mathfrak{S}}(B)$ is not sylleptic in general. Nonetheless, under favourable circumstances, this is in fact the case. We begin by recalling the following definition from section 5.3 of \cite{Cr}.

\begin{definition}\label{def:syllepticcenter}
Let $\mathfrak{S}$ be a sylleptic fusion 2-category. The symmetric center of $\mathfrak{S}$, denoted by $\mathcal{Z}_{(3)}(\mathfrak{S})$ is the full sub-2-category of $\mathfrak{S}$ on those objects $C$ such that $$\sigma_{D,C}\circ b_{C,D} = b_{C,D}\circ \sigma_{C,D}$$ for every $D$ in $\mathfrak{S}$.
\end{definition}

\begin{rem}
It follows immediately from the definitions that $\mathcal{Z}_{(3)}(\mathfrak{S})$ is a (semi-strict) symmetric monoidal 2-category (see also theorem 5.2 of \cite{Cr}).
\end{rem}

\begin{proposition}\label{prop:SymmetricAlgebraCenterModule}
Let $\mathfrak{S}$ be a Karoubi complete sylleptic monoidal 2-category, and $B$ a symmetric separable algebra in $\mathcal{Z}_{(3)}(\mathfrak{S})$. Then, $\mathbf{Mod}_{\mathfrak{S}}(B)$ is a sylleptic monoidal 2-category.
\end{proposition}

\begin{proof}
Without loss of generality, we may assume that $\mathfrak{S}$ is semi-strict. We have already endowed the 2-category $\mathbf{Mod}_{\mathfrak{S}}(B)$ with a braided monoidal structure. Moreover, using the notation of the proof of theorem \ref{thm:SymmetricAlgebraModule}, for every right $B$-modules $M$ and $N$ in $\mathfrak{S}$, we can consider the following right $B$-module 2-isomorphism $$\begin{tikzcd}
M\Box N \arrow[dd, equal] \arrow[rd, "b"] \arrow[r, "t"] & M\Box_B N \arrow[rd, "\widetilde{b}"]                    &                                       \\
{} \arrow[r, "\sigma", Rightarrow, shorten > = 2ex, shorten < = 3ex]                                  & N\Box M \arrow[ld, "b"'] \arrow[r, "t"] \arrow[u, "\xi", Rightarrow] & N\Box_B M \arrow[ld, "\widetilde{b}"] \\
M\Box N \arrow[r, "t"'] \arrow[rru, "\xi", Rightarrow, shorten > = 8ex, shorten < = 8ex]              & M\Box_B N.                                                &                                      
\end{tikzcd}$$ The above right $B$-module 2-isomorphism is $B$-balanced. In order to see this, we use the diagrams depicted in appendix \ref{subsub:proofSymmetricAlgebraCenterModule}. Figure \ref{fig:sigmaxixibalanced1} depicts the left hand-side of equation (\ref{eqn:balance2morphism}) of definition \ref{def:balanced2morphism} for the above 2-isomorphism. Applying equation (\ref{eqn:balance2morphism}) for $\xi$ on the blue coupons brings us to figure \ref{fig:sigmaxixibalanced2}. By inserting the definition of $\tau^{t\circ b}$ given in the proof of theorem \ref{thm:SymmetricAlgebraModule}, we arrive at figure \ref{fig:sigmaxixibalanced3}. Then, using equation (\ref{eqn:balance2morphism}) for $\xi$ on the blue coupons leads us to figure \ref{fig:sigmaxixibalanced4}. Inserting the definition of $\tau^{t\circ b}$ once again, we get to figure \ref{fig:sigmaxixibalanced5}. In order to get to figure \ref{fig:sigmaxixibalanced6}, we first use the equation given in definition \ref{def:syllepticcenter} on the blue coupons and the strand immediately on top of it, and then move the left most coupon labeled $\sigma$ along the green arrow. Then, applying equation (\ref{eqn:syllepsisaxiom1}) on the blue coupons brings to figure \ref{fig:sigmaxixibalanced7}. Using equation (\ref{eqn:syllepsisaxiom2}) on the blue coupons, followed by moving the freshly created coupon labeled $\sigma$ down along the green arrow, and cancelling the pair of red coupons brings us to figure \ref{fig:sigmaxixibalanced8}. But, figure \ref{fig:sigmaxixibalanced8} depicts the right hand-side of equation (\ref{eqn:balance2morphism}), so the proof of the claim is finished.

Then, thanks to the 2-universal property of the relative tensor product, this yields a 2-isomorphism $\widetilde{\sigma}_{M,N}$ as in the diagram below $$\begin{tikzcd}[sep=small]
M\Box_BN \arrow[rrrr, equal] \arrow[rrdd, "\widetilde{b}"'] &  & {} \arrow[dd, Rightarrow, "\widetilde{\sigma}", near start, shorten > = 1ex] &  & M\Box_BN. \\
                                   &  &                           &  &   \\
                                   &  & N\Box_B M \arrow[rruu, "\widetilde{b}"']     &  &  
\end{tikzcd}$$ Further, it follows from the 2-universal property of the realtive tensor product that the collection of the 2-isomorphisms $\widetilde{\sigma}_{M,N}$ for varying $M$ and $N$ defines an invertible modification. Finally, one has to check that $\widetilde{\sigma}$ defines a syllepsis on the braided monoidal 2-category $\mathbf{Mod}_{\mathfrak{S}}(B)$, i.e. that the equations given in figure C.15 and C.16 of \cite{SP} hold. This follows from the 2-universal property of the relative tensor product over one and two algebras explained in section 3 of \cite{D8}.
\end{proof}

We now consider the case when $\mathfrak{S}$ is symmetric monoidal.

\begin{corollary}
Let $\mathfrak{S}$ be a Karoubi complete symmetric monoidal 2-category, and $B$ a symmetric separable algebra in $\mathfrak{S}$. Then, $\mathbf{Mod}_{\mathfrak{S}}(B)$ is a symmetric monoidal 2-category.
\end{corollary}
\begin{proof}
If $\mathfrak{S}$ is symmetric, then $\mathcal{Z}_{(3)}(\mathfrak{S})=\mathfrak{S}$, which implies that $\mathbf{Mod}_{\mathfrak{S}}(B)$ is sylleptic. Further, it follows from the definition of the syllepsis, that $\mathbf{Mod}_{\mathfrak{S}}(B)$ is in fact symmetric if $\mathfrak{S}$ is symmetric.
\end{proof}

\begin{lemma}\label{lem:freemodulesylleptic}
Let $\mathfrak{S}$ be a Karoubi complete sylleptic monoidal 2-category, and $B$ a symmetric separable algebra in $\mathcal{Z}_{(3)}(\mathfrak{S})$. Then, the braided monoidal functor $F:\mathfrak{S}\rightarrow \mathbf{Mod}_{\mathfrak{S}}(B)$ of proposition \ref{lem:freemodulemonoidal} is sylleptic. In particular, if $\mathfrak{S}$ is symmetric, then $F$ is symmetric.
\end{lemma}
\begin{proof}
The first part follows from the construction and the 2-universal property of the relative tensor product over $B$. The last part is immediate as a symmetric monoidal 2-functor is nothing but a sylleptic monodial 2-functor between symmetric monoidal 2-categories (see definition 2.5 of \cite{SP}).
\end{proof}

\begin{rem}
Analogously to what was noted in remark \ref{rem:generaltensor}, the proofs of all the above results in this section only used the existence of the relative tensor product over $B$ for any $B$-modules in $\mathfrak{S}$.
\end{rem}

Finally, if $\mathfrak{S}$ is a sylleptic multifusion 2-category, proposition \ref{prop:SymmetricAlgebraCenterModule} can be strengthened. We begin by the following lemma.

\begin{lemma}\label{lem:symmetriccenter}
Let $\fS$ be a sylleptic fusion 2-category. Then, its symmetric center $\mathcal{Z}_{(3)}(\fS)$ is generated under direct sums by the union of some of the connected components of $\fS$. In particular, it is a symmetric fusion 2-category, and it contains the connected components of the identity of $\fS$.
\end{lemma}
\begin{proof}
Observe that, by definition, $\mathcal{Z}_{(3)}(\fS)$ is a full sub-2-category of $\fS$. Further, note that $\mathcal{Z}_{(3)}(\fS)$ is closed under taking direct sums. Now, let $S$ be an object of $\mathcal{Z}_{(3)}(\fS)$. We wish to prove that if $T$ is a simple object of $\mathfrak{S}$ given by the splitting of a 2-condensation monad on $S$, then $T$ is in $\mathcal{Z}_{(3)}(\fS)$. Given an arbitrary object $C$ in $\mathfrak{S}$, it follows from the 2-universal property of the splitting of a 2-condensation monad that the syllepsis $\sigma_{T,C}$ and $\sigma_{C,T}$ are completely determined by $\sigma_{S,C}$ and $\sigma_{C,S}$. But, by hypothesis, we have $\sigma_{C,S}\circ b_{S,C} = b_{S,C}\circ \sigma_{C,S}$, so that $\sigma_{C,T}\circ b_{T,C} = b_{T,C}\circ \sigma_{C,T}$, which proves the claim. The second part follows from the observation that a connected component of a finite semisimple 2-category is necessarily a finite semisimple 2-category.
\end{proof}

\begin{proposition}
Let $\mathfrak{S}$ be a sylleptic multifusion 2-category, and $B$ a symmetric separable algebra in $\mathcal{Z}_{(3)}(\mathfrak{S})$. Then, we have $$\mathcal{Z}_{(3)}(\mathbf{Mod}_{\mathfrak{S}}(B))\simeq \mathbf{Mod}_{\mathcal{Z}_{(3)}(\mathfrak{S})}(B).$$
\end{proposition}
\begin{proof}
It follows from the construction that the syllepsis on $\mathbf{Mod}_{\mathfrak{S}}(B)$ is constructed from the syllepsis on $\mathfrak{S}$. In particular there is a symmetric monoidal inclusion $\mathbf{Mod}_{\mathcal{Z}_{(3)}(\mathfrak{S})}(B)\subseteq \mathcal{Z}_{(3)}(\mathbf{Mod}_{\mathfrak{S}}(B)).$ On the other hand, the free 2-functor $F:\fS\rightarrow \mathbf{Mod}_{\mathfrak{S}}(B)$ is sylleptic monoidal. In particular, for any object $C$ of $\fS$, $F(C)$ is in $\mathcal{Z}_{(3)}(\mathbf{Mod}_{\mathfrak{S}}(B))$ if and only if $C$ is in $\mathcal{Z}_{(3)}(\mathfrak{S})$. But, every object of $\mathbf{Mod}_{\mathfrak{S}}(B)$ is the splitting of a 2-condensation monad supported on an object of the form $F(C)$ for some $C$ in $\fS$ by lemma 3.1.1 of \cite{D7}. Further, $\mathcal{Z}_{(3)}(\mathbf{Mod}_{\mathfrak{S}}(B))$ is a union of connected components of $\mathbf{Mod}_{\mathfrak{S}}(B)$ by lemma\ref{lem:symmetriccenter} above, so that $\mathcal{Z}_{(3)}(\mathbf{Mod}_{\mathfrak{S}}(B))\simeq \mathbf{Mod}_{\mathcal{Z}_{(3)}(\mathfrak{S})}(B)$ as desired.
\end{proof}

\section{Specific 2-Category of Modules}\label{section:ExampleCats}

In this section, we will examine the 2-categories of right modules associated to specific algebras. This can be thought of as condensing a 3-condensation monad. In order to be applicable to physical theories, we will consider the cases when the ambient 2-category is either totally disconnected or connected. In both cases, we will work bosonically and
fermionically, where the later means that we work with super 2-categories. A subset of surface operators can be assembled to form a separable algebra as in the previous section; we may thus apply the theorems above to understand the effect of the condensation. Throughout, we work over the complex numbers (or any algebraically closed field of characteristic zero), we use $G$ to denote a finite group, and $E$ to denote a finite abelian group. 

\subsection{Totally Disconnected 2-Category}

\subsubsection{Bosonic case}Starting with the simplest case, suppose that the fusion 2-category of surface operators and their interactions is given by $2\Vec[G]$, the 2-category of $G$-graded 2-vector spaces. In particular, the (equivalence classes of) simple objects are given by $\mathbf{Vect}_g$ with $g\in G$. We can consider the algebra $\Vec{[G]}$ in $2\Vec[G]$
given by $\boxplus_{g\in G}\mathbf{Vect}_g$, the sum of the equivalence classes of simple objects.

\begin{lemma}
The left $2\Vec[G]$-module 2-category $2\Vec$ is equivalent to $\Mod_{{2\Vec[G]}}(\Vec{[G]})$, where $\Vec{[G]}$ is the fusion 1-category of $G$-graded vector spaces viewed as an algebra in $2\Vec[G]$ with the canonical grading. Further, $\Vec{[G]}$ is a separable algebra.
\end{lemma}

\begin{proof}
It is easy to check directly that $\Mod_{{2\Vec[G]}}(\Vec{[G]})\simeq 2\Vec$ as left $2\Vec[G]$-module 2-categories. Further, one can check directly that $\Vec{[G]}$ is a separable algebra in $2\Vec[G]$. Alternatively, this follows from theorem 3.2.4 and corollary 3.3.7 of \cite{D7}.
\end{proof}


Before moving on to the general case, we establish the following technical result. Recall from \cite{D4} that a left $2\Vec[G]$-module 2-category is 2-category equipped with a left action by $2\Vec[G]$. Note that this is equivalent to the data of an action of the group $G$. In particular, given a left $2\Vec[G]$-module 2-category $\mathfrak{M}$, we can consider the 2-category $\mathbf{LMod}_{\mathfrak{M}}(\Vec[G])$ of left $\Vec[G]$-modules in $\mathfrak{M}$, given by gauging the $G$-action on $\mathfrak{M}$. If $\mathfrak{M}$ is a finite semisimple 2-category, the $G$-action permutes the set of connected components of $\mathfrak{M}$.

\begin{proposition}\label{prop:equivariantization}
Let $\mathfrak{M}$ be a finite semisimple left $2\Vec[G]$-module 2-category. Then, we have $$\pi_0(\mathbf{LMod}_{\mathfrak{M}}(\Vec[G])) \cong \pi_0(\mathfrak{M})/G.$$
\end{proposition}
\begin{proof}
We claim that it suffices to prove this result for $\mathfrak{M}$ an indecomposable finite semisimple left $2\Vec[G]$-module 2-category. Namely, it follows from lemma 5.2.3 of \cite{D8} that every finite semisimple left $2\Vec[G]$-module 2-category $\mathfrak{M}$ can be decomposed into a finite direct sum $\mathfrak{M}\simeq \boxplus_{i=1}^n\mathfrak{M}_i$ of indecomposable ones. From this, it follows that there is a bijection $\pi_0(\mathfrak{M})\cong \coprod_{i=1}^n\pi_0(\mathfrak{M}_i)$ of sets compatible with the $G$-actions. This establishes the claim of sufficiency.

Now, note that it follows from the definition that a finite semisimple left $2\Vec[G]$-module 2-category is indecomposable if and only if the action of $G$ on $\pi_0(\mathfrak{M})$ is transitive. Thus, it only remains to prove that if $\mathfrak{M}$ is an indecomposable finite semisimple left $2\Vec[G]$-module 2-category, then $\pi_0(\mathbf{Mod}_{\mathfrak{M}}(\Vec[G]))=*$.

To see this, note that thanks to theorem 5.1.2 of \cite{D4}, there exists an algebra $A$ in $2\Vec[G]$ such that $\mathfrak{M}\simeq \mathbf{Mod}_{2\Vec[G]}(A)$. Furthermore, by theorem 5.4.7 of \cite{D4}, the algebra $A$ is in fact rigid. But rigid algebras in $2\Vec[G]$ are precisely $G$-graded multifusion 1-categories, so that $A$ is an $G$-graded multifusion 1-category. Moreover, as $\mathfrak{M}$ is indecomposable, $A$ is indecomposable as a $G$-graded multifusion 1-category (see corollary 5.2.7 of \cite{D8}).

By inspection, there are equivalences of 2-categories $$\mathbf{LMod}_{\mathfrak{M}}(\Vec[G])\simeq \mathbf{Bimod}_{2\Vec[G]}(\Vec[G],A)\simeq \mathbf{Mod}_{2\Vec}(A),$$ where, on the right hand-side, we view $A$ as a multifusion 1-category. Thus, by proposition 2.3.5 of \cite{D1}, it is enough to prove that $A$ is indecomposable as a multifusion 1-category. (A multifusion 1-category is ``connected'' in the sense of definition 2.3.1 of \cite{D1} if and only if it is indecomposable.) Finally, observe that a decomposition of $A$ into a direct sum of two non-zero multifusion 1-categories would automatically be compatible with the $G$-grading. This is impossible by construction so we are done.
\end{proof}

If $G=E$ is a finite abelian group, then $2\Vec[E]$ is braided fusion 2-category. Further, the algebra $\Vec{[E]}$ is actually braided. It is therefore sensible to consider the case when the 2-category of all surfaces is a braided fusion 2-category $\mathfrak{B}$, equipped with a braided monoidal inclusion $2\Vec[E]\subseteq \mathfrak{B}$. This allows us to view the separable algebra $\Vec[E]$ in $2\Vec[E]$ as living in $\mathfrak{B}$, and we can investigate the properties of 2-category obtained by the condensation of $\Vec[E]$ in $\mathfrak{B}$. The following result follows from theorem \ref{thm:multifusion2cat}, the above proposition, and lemma \ref{lem:freemodulemonoidal}.

\begin{corollary}\label{cor:disconnectedbraided}
Given $\fB$ a braided fusion 2-category and $2\Vec{[E]} \subseteq\fB$ a braided monoidal inclusion, the 2-category $\Mod_{\fB}(\Vec{[E]})$ obtained by condensing $\Vec[E]$ is a fusion 2-category with $\pi_0(\Mod_{\fB}(\Vec{[E]})) \cong \pi_0(\fB)/E$. Moreover, the canonical 2-functor $\fB\rightarrow \Mod_{\fB}(\Vec{[E]})$ is monoidal.
\end{corollary}

In particular, the condensation reorganizes the 2-category $\fB$ by identifying the connected components of surfaces which are related by the
action of $E$. This is effectively gauging the $E$ action on the components. The resulting fusion 2-category is in general not braided.

\begin{example}
Consider $2\Vec[\bZ_4]$, with simple objects labeled by $\{\mathbf{Vect}_0,\mathbf{Vect}_1,\mathbf{Vect}_2,\mathbf{Vect}_3\}$ and fusion given by addition mod 4. Suppose we condense the algebra $\mathbf{Vect}_0\boxplus\mathbf{Vect}_2$, which is $\Vec[\bZ_2]$, the simple modules are then given by $\mathbf{Vect}_0\boxplus\mathbf{Vect}_2$, and $\mathbf{Vect}_1\boxplus\mathbf{Vect}_3$. As there is no 1-morphism between them, $\Mod_{2\Vec[\bZ_4]}(\Vec{[\bZ_2]})$ has two connected components. On the other hand, one sees that $\pi_0(2\Vec[\bZ_4])/\bZ_2$ has the same two connected components.
\end{example}

\begin{rem}
We give an example for which the 2-functor in corollary \ref{cor:disconnectedbraided} is not necessarily braided, take $\mathfrak{B}=\mathcal{Z}(2\Vec[\bZ_2])$, the Drinfeld center of $2\Vec[\bZ_2]$, equipped with the canonical inclusion $2\Vec [\mathbb{Z}_2]\subseteq \mathcal{Z}(\mathbf{2Rep}(\mathbb{Z}_2))$. We can then condense the algebra $\Vec[\bZ_2]$, and get $$\mathbf{Mod}_{\mathcal{Z}(2\Vec[\bZ_2])}(\Vec [\mathbb{Z}_2])\simeq \mathbf{2Rep}(\mathbb{Z}_2).$$ Further, the monoidal 2-functor $\mathcal{Z}(2\Vec[\bZ_2])\rightarrow\mathbf{Mod}_{\mathcal{Z}(\mathbf{2Rep}(\mathbb{Z}_2))}(\Vec [\mathbb{Z}_2])$ of lemma \ref{lem:freemodulemonoidal} is identified with the monoidal forgetful  2-functor $\mathcal{Z}(\mathbf{2Rep}(\mathbb{Z}_2))\rightarrow\mathbf{2Rep}(\mathbb{Z}_2)$, which is not braided.
\end{rem}

The next result follows from proposition \ref{prop:SymmetricAlgebraCenterModule}, lemma \ref{lem:freemodulesylleptic}, and proposition \ref{prop:equivariantization}.

\begin{corollary}
Let $\fS$ be a sylleptic fusion 2-category, with an inclusion $2\Vec{[E]}\subseteq \mathcal{Z}_{(3)}(\fS)$, then $\Mod_{\fS}(\Vec{[E]})$ is a sylleptic fusion 2-category such that $\pi_0(\Mod_{\fS}(\Vec{[E]})) \cong \pi_0(\fS)/E$. Furthermore, the canonical monoidal 2-functor $\fS\rightarrow \Mod_{\fS}(\Vec{[E]})$ is sylleptic.
\end{corollary}

\subsubsection{Fermionic Case} We mirror the bosonic case and first consider the fusion 2-category
$2\SVec[G]$ of $G$-graded super 2-vector spaces. In order to condense $2\SVec[G]$ to $2\SVec$, it is enough to consider the bosonic algebra $\Vec[G]$ given by the canonical monoidal inclusion $2\Vec[G]\subseteq 2\SVec[G]$. By direct inspection, we find that $\Mod_{2\SVec[G]}(\Vec[G])\simeq 2 \SVec$.

Let us now comment on the braided case. Namely, if $G=E$ is a finite abelian group, then $2\SVec[E]$ is a braided fusion 2-category. We can therefore consider $\mathfrak{B}$ a braided fusion 2-category containing $2\SVec[E]$. But, the inclusion $2\Vec[E]\subseteq 2\SVec[E]$ is braided, so this is exactly in the setup of corollary \ref{cor:disconnectedbraided}. Similar remarks holds for the sylleptic and symmetric cases.

\subsection{Connected Category}\label{subsection:Connected}
Let $\fB$ be a braided fusion 2-category, then $\End_{\fB}(\mathbb{1})$, the endomorphisms of the identity surface, is a symmetric fusion 1-category, so that $\fB^0 = \mathbf{Mod}(\End_{\fB}(\mathbb{1}))$ is a symmetric fusion 2-category (see \cite{D2}). Here, $\mathfrak{B}^0$ denotes the identity component and is a prime candidate for a condensation. 

\subsubsection{Bosonic Case}\label{subsub:Connectedbosonic} Suppose that $\mathfrak{B}^0=2\rRep(G)$, i.e. the surfaces in the identity component of $\mathfrak{B}$ form the fusion 2-category $2\rRep(G)$. Here we think of $2\rRep(G)$ as the 2-category of finite semisimple 1-categories equipped with a $G$-action. One such object is given by $\Vec[G]$ with the canonical $G$-action. In this description, the monoidal product of two finite semisimple 1-categories $\mathcal{C}$ and $\mathcal{D}$ equipped with $G$-actions is given by their Deligne tensor product $\mathcal{C} \boxtimes \mathcal{D}$ equipped with the diagonal $G$-action. The fusion 2-category $2\rRep(G)$ is connected, which means that all the surfaces arise as networks of lines. We write $\varphi$ for the symmetric algebra $\mathbf{Fun}(G,\Vec)$ in $2\rRep(G)$. We note that the underlying object of $\varphi$ is $\Vec[G]$. In the setting of fusion 1-categories, this corresponds to considering the symmetric algebra $\mathbb{C}[G]^*$ inside $\rRep(G)$. A module for $\varphi$ is thus a way for the lines to end at the boundary.

We also point out that there is another model for $2\rRep(G)$, given by $\Mod(\rRep(G))$ (see lemma 1.3.8 of \cite{D8}). In the fermionic case, only this second model is available. It is therefore necessary to give an alternative description of $\varphi$ in this model. The symmetric fusion 1-category $\Vec$ equipped with the canonical symmetric monoidal functor $\rRep(G)\rightarrow \Vec$ defines a symmetric algebra in $\Mod(\rRep(G))$. This algebra is separable thanks to theorem 3.2.4 and proposition 3.3.3 of \cite{D7} and theorem 2.3 of \cite{ENO1}. Moreover, under the equivalence of lemma 1.3.8 of \cite{D8}, the algebra $\Vec$ in $\Mod(\rRep(G))$ corresponds to the algebra $\varphi$ in the first model. It follows that $\Mod_{2\rRep(G)}(\varphi) \simeq \Mod_{\Mod(\rRep(G))}(\Vec)\simeq 2\Vec$.

\begin{proposition}\label{prop:connectedSF}
Let $\mathfrak{B}$ be a braided fusion 2-category with $2\rRep(G)\simeq\fB^0$ as braided fusion 2-categories. Then, condensing the braided separable algebra $\varphi$ in $\mathfrak{B}$ yields a strongly fusion 2-category $\Mod_{\fB}(\varphi)$ equipped with a monoidal 2-functor $\fB\rightarrow \Mod_{\fB}(\varphi)$.
\end{proposition}
\begin{proof}
All but the strongly fusion part follow from theorem \ref{thm:multifusion2cat}. We claim that $\Mod_{\fB}(\varphi)^0=2\Vec$, so that $\Mod_{\fB}(\varphi)$ is strongly fusion. Note that $\varphi$ is an algebra in $2\rRep(G)\simeq\fB^0$. By corollary 2.3.6 of \cite{D2}, this implies that the underlying object in $\fB$ of any simple right $\varphi$-module is supported in a single connected component of $\fB$. This shows that $\Mod_{\fB}(\varphi)^0\simeq\Mod_{2\rRep(G)}(\varphi)^0\simeq 2\Vec$. This finishes the proof of the claim.
\end{proof}

\begin{rem}
In the fusion 2-category $2\rRep(G)$ the algebra $\varphi$ is actually symmetric, but we can not view $\varphi$ as a symmetric algebra in $\mathfrak{B}$; this requires extra data in the ambient braided fusion 2-category $\mathfrak{B}$. Therefore $\varphi$ is treated as a braided algebra when considered in $\mathfrak{B}$.
\end{rem}

We give a physical explanation as to why
condensing in the identity component in proposition \ref{prop:connectedSF} was sufficient to make $\mathfrak{B}$ strongly fusion: the objects in the identity component of $\mathfrak{B}$ are related to the identity surface by 2-condensations but if the identity component was condensed to just 2$\Vec$ via $\varphi$, then all the 1-morphisms are trivial, hence the 2-category is strongly fusion.

\begin{rem}
Categorifying the main result of \cite{Ki} and \cite{Mu}, we expect that if $\mathfrak{B}$ is braided fusion 2-category with $2\rRep(G)\simeq\fB^0$, then the fusion 2-category $\mathbf{Mod}_{\mathfrak{B}}(\varphi)$ admits a $G$-crossed braided structure.
\end{rem}

Let $\mathfrak{S}$ be a sylleptic multifusion 2-cateogry. As a consequence of lemma \ref{lem:symmetriccenter}, we find that any inclusion $2\rRep(G) \subseteq \mathfrak{S}$ of sylleptic fusion 2-categories automatically includes in the symmetric center of $\mathfrak{S}$. Namely, $2\rRep(G)$ is necessarily contained in the component of the identity of $\mathfrak{S}$. Combing this observation with proposition \ref{prop:SymmetricAlgebraCenterModule} and lemma \ref{lem:freemodulesylleptic} yields the following result.

\begin{corollary}\label{cor:stronglysylleptic}
Let $\fS$ be a sylleptic fusion 2-category. Suppose that there is an inclusion $\fS\simeq2\rRep(G)$, then $\Mod_{\fS}(\varphi)$ is a sylleptic strongly fusion 2-category. Furthermore, the canonical monoidal 2-functor $\fS\rightarrow \Mod_{\fS}(\varphi)$ is sylleptic.
\end{corollary}

\begin{rem}
We make a small physical point regarding the above corollary. Consider a setting in  (3+1)d but \textit{not} limited to considering only topological theories. The surface operators can be nontrivial even if the line operators have been condensed, as in the situation of corollary \ref{cor:stronglysylleptic}. If we are in a purely topological (3+1)d setting, then the there are actually no surface operators either because surfaces detect lines in this dimension. This means we are just in a situation of bosonic Dijkgraaf-Witten theory.
\end{rem}

\subsubsection{Fermionic Case} We consider the fusion 2-category $2\rRep(G,z):=\Mod(\rRep(G,z))$, where $z$ is an emergent fermion in $G$, that is a central element of order 2. We are viewing $2\rRep(G,z)$ as so because there is no fermionic analogue of the model for $2\rRep(G)$ that was used in \S\ref{subsub:Connectedbosonic}. We define the symmetric separable algebra $\varphi := \SVec$ in $2\rRep(G,z)$. More precisely, $\varphi$ denotes $\SVec$ equipped with the canonical forgetful symmetric monoidal functor $\rRep(G,z)\rightarrow \SVec$. Let us examine the result of condensing $\varphi$. In this case, there is no obstruction to condensing to the vacuum.

\begin{lemma}
As left $2\rRep{(G,z)}$-module 2-categories, we have $2\SVec \simeq \mathbf{Mod}_{2\rRep{(G,z)}}(\varphi)$, where $\SVec$ is viewed as an algebra in $2\rRep{(G,z)}$ via $\rRep{(G,z)}\rightarrow \SVec$.
\end{lemma}
\begin{proof}
This follows from example 3.2.5 of \cite{D4}.
\end{proof}

Physically, we find that condensing $\varphi$ gives a local fermion. The next proposition follows using a variant of the proof of proposition \ref{prop:connectedSF}, with a slight change to $\varphi$.

\begin{proposition}
Let $\mathfrak{B}$ be a braided fusion 2-category, and assume $\fB^0\simeq2\rRep(G,z)$ as braided fusion 2-categories. Then, condensing the algebra $\varphi = \SVec$ in $\mathfrak{B}$ yields a fermionic strongly fusion 2-category $\Mod_{\fB}(\varphi)$ equipped with a monoidal 2-functor $\fB\rightarrow \Mod_{\fB}(\varphi)$.
\end{proposition}

\section{Strongly Fusion Computations}\label{section:StronglyFusion}

\subsection{Braided Strongly Fusion 2-Categories}
In the previous section, we have seen examples of strongly fusion 2-categories arising from condensations. It was shown in \cite{JFY} that such 2-categories have grouplike fusion rules. Said differently, a strongly fusion 2-category is a ``grouplike'' extension of operators in different dimensions \cite{JF2,JFY2}. In particular, their classification essentially boils down to a cohomology computation problem. 
We now consider the case where our fusion 2-category is strongly fusion and only braided. For instance, this is what happens to a sylleptic fusion 2-category when we condense the algebra $\varphi$ in the subcategory $\mathfrak{B}^0=2\rRep(G)$. Fermionic braided strongly fusion 2-categories are classified by supercohomology \cite{Wang:2017moj} and we expect that the cases we discuss here cover all the examples of braided strongly fusion 2-categories, but we do not prove this fact. Namely, in general, one ought to consider supercohomology with twisted coefficients, but we expect that this is not necessary for braided fermionic strongly fusion 2-categories. On the other hand, braided bosonic strongly fusion 2-categories with finite abelian group of surfaces given by $E$ are completely classified by $\rH^5(E[2];\mathbb{C}^\times)$. This holds because this cohomology theory has no twisted variant.

The classification of the physical theories described by braided strongly fusion 2-categories proceeds by identifying those fusion 2-categories that are related by a topological boundary. More precisely, fixing a finite abelian group of surfaces $E$, the associated physical theories are classified by generalized cohomology. In the fermionic case, the relevant spectrum of coefficients is $\mathcal{SW}^\bullet(\pt)$, the super--Witt spectrum \cite{JF}. Its homotopy groups in low degrees are recalled below in \eqref{eq:SWpt}. In the bosonic case, the classification requires twisted equivariant cohomology. We now discuss these computations in more detail.

\subsubsection{Fermionic Case}

Let $\fB$ be a braided fermionic strongly fusion 2-category, and write $E$ for the finite abelian group of connected components. Physically, $E$ is the group of ``fundamental" surfaces in $\mathfrak{B}$ that do not arise as condensations. Further, such braided fusion 2-categories can be constructed by deforming the coherence structure of $2\SVec{[E]}$ using a class in the super-cohomology group $\SH^5(E[2])$. Here and in what follows, $E[2]$ denotes the second Eilenberg-MacLane space of $E$, and we note that the number in brackets denotes the codimension associated to the objects with fusion rules given by the group $E$. In fact, all braided fermionic strongly fusion 2-categories arise via this construction, but we do not prove this fact. On the other hand, the (3+1)d theory associated to $\fB$ has no codimension one operators that do not arise through a condensation. By remote detectability \cite{JF}, which says that every object must link topologically with another object of the appropriate dimension, this is the same as assuming that there are no nontrivial point operators in the theory. Then, the obstruction to condensing the theory associated to $\fB$ to the vacuum is given by a class in $\mathcal{SW}^\bullet(E[2])$. Now, if the group $\mathcal{SW}^\bullet(E[2])$ vanishes the theory associated to $\fB$ is automatically Morita equivalent to the vacuum. Our goal is to understand for which abelian groups $E$ the cohomology group $\mathcal{SW}^5(E[2])$ does not vanish. More precisely, there is a canonical map $\SH^5(E[2])\rightarrow \mathcal{SW}^\bullet(E[2])$, which corresponds to taking the theory associated to a braided fermionic strongly fusion 2-category. We argue that the image of this map is non-trivial in general.

Since the fusion 2-category $\fB$ is strongly fusion, 
there are no nontrivial lines, but we still have $\{1,f\}$ in $\SVec$ from the fermionic nature of the 2-category. We denote the condensation surface arising from $f$ as $c$ which has fusion rule $c^2 = \mathbb{1}$.
The content of the fusion 2-category forms a higher group extension
\begin{align}
        (\bC ^\times [4] .\underset{\{1,f\}}{\underbrace{\bZ_2}}[3] .\underset{\{\mathbb{1},c\}}{\underbrace{\bZ_2}}[2]).E[2]\,,
\end{align}
where the component $\bC^\times [4]$ means ``three form $\bC^\times$ symmetry". Such extensions are classified by $\SH^5(E[2])$, which can be computed with the knowledge that the supercohomology of a point is built out of three layers: 
\begin{equation}
      \SH^0(\pt)= \mathbb{C}^\times\,,\quad  \SH^1(\pt)= \bZ_2\,, \quad  \SH^2(\pt)= \bZ_2\,,
\end{equation}
We note in passing that these groups agree with the first three layers of spin cobordism. Then, there is a canonical map $\SH^{\bullet}\rightarrow \mathcal{SW}^{\bullet}$. Assuming that $\mathfrak{B}$ is classified by a class in $\SH^5(E[2])$, the associated fermionic theory can be condensed to the vacuum exactly if the image of this class in $\mathcal{SW}^\bullet(E[2])$ is trivial. In order to understand for which groups $E$ this can happen, we use the following Atiyah-Hirzebruch spectral sequence
\begin{equation}\label{eq:SWA2}
    \rH^i(E[2];\, \mathcal{SW}^j(\pt)) \Rightarrow  \mathcal{SW}^{i+j}(E[2])\,.
\end{equation}
The homotopy groups of $\SW^\bullet(\pt)$ in low degrees are given by 
\begin{align}\label{eq:SWpt}
    \mathcal{SW}^0(\pt) &=\bC^\times \,,\quad  \mathcal{SW}^1(\pt) =\bZ_2 \,, \quad \SW^2(\pt) =\bZ_2\,, \\\notag 
    \mathcal{SW}^3(\pt) &=0\,,\quad \mathcal{SW}^4(\pt) =\mathcal{SW}\,,\quad \mathcal{SW}^5(\pt) =0\,,\quad \mathcal{SW}^6 (\pt)=0\,.
\end{align}
In degree 4, $\mathcal{SW}$ gives the Witt group of slightly degenerate braided fusion 1-categories. If $E$ has no 2-torsion, then we find that $\mathcal{SW}^5(E[2])=\SH^5(E[2])=\rH^5(E[2];\bC^\times)$. But, it follows from \cite{EilenbergMacLane} that the right most group is trivial, so that there are non-trivial theories in this case. On the other hand, we can assume that $E$ is 2-torsion. Then, we see that in total degree 5, there are interesting non-zero contributions to the $E_2$-page of the spectral sequence. We first consider $E=\bZ_{2^k}$.
The $E_2$ page for \eqref{eq:SWA2} is then given by 
\begin{equation}\label{E2pageK2}
    E^{ij}_2 = \,\begin{array}{c|ccccccccc}
      j\\
      \\
      0&  0  &0&\ldots \\
    0&  0  &0&\ldots \\
     \mathcal{SW} &   \mathcal{SW}& 0& \hom(\mathcal{SW},\bZ_2)&\ldots \\
     0& 0&  0 & 0 & 0 & \ldots\\
     \bZ_2&\bZ_2&0&\hom(E,\bZ_2)&\Ext(E,\bZ_2)&\Quad{(E,\bZ_2)} &\ldots \\
     \bZ_2&\bZ_2&0&\hom(E,\bZ_2)&\Ext(E,\bZ_2)&\Quad{(E,\bZ_2)}&\Ext(E,\hom(E,\bZ_2))&  \\
     \bC^\times & \bC^\times& 0 & \widehat{E} & 0 & \Quad{(E,\bC^\times)}& \Ext(E,\hom(E,\bC^\times))   \\
     \hline 
     & 0 & 1 & 2 & 3 & 4 & 5 &  \quad i\,,
\end{array}
\end{equation}
where $\Quad$ denotes the group of quadratic forms. In addition, the $d_2$ differentials are given by
 \begin{align}\label{eq:d2differentialuntwisted}
     d_2:& E^{i,2}_2=\rH^i(\bZ_2[2]\,; \bZ_2) \to E^{i+2,1}_2=\rH^{i+2}(\bZ_2[2]\,;\bZ_2) \quad&& X \mapsto \Sq^2 X \\\notag 
     d_2:& E^{i,1}_2=\rH^i(\bZ_2[2]\,; \bZ_2) \to E^{i+2,0}_2=\rH^{i+2}(\bZ_2[2]\,; \bC^\times) \quad&& X \mapsto (-1)^{\Sq^2 X} \,.
\end{align}
This implies that the $E_3$ page is given by 
\begin{equation}
    E^{ij}_3 = \,\begin{array}{c|cccccccccc}
      j\\
      \\
      0&  0  &0&\ldots \\
    0&  0  &0&\ldots \\
     \mathcal{SW} &    \mathcal{SW} & 0&  \mathcal{SW}_2&\ldots \\
     0& 0&  0 & 0 & 0 &  \ldots\\
     \bZ_2&\bZ_2&0&0&0& \ldots \\
     \bZ_2&\bZ_2&0&0&0&0&\ldots \\
     \bC^\times & \bC^\times& 0 & \bZ_2 & 0 & \bZ_2& 0& \bZ_2 &  \\
     \hline 
     & 0 & 1 & 2 & 3 & 4 & 5 & 6&\quad i\,.
\end{array}
\end{equation}
Therefore, $\mathcal{SW}^5(\bZ_{2^k}[2])=0$, so that the theory associated to $\fB$ with $E=\bZ_{2^k}$ can be condensed to the vacuum.

If $E$ is a product of groups, we use the fact that for any generalized cohomology theory $h^\bullet$ computed on pointed spaces $X$ and $Y$, we have
\begin{equation}\label{eq:prodgroup}
    h^\bullet(X\times Y) = h^\bullet(\pt) \oplus \widetilde{h}^\bullet(X) \oplus \widetilde{h}^\bullet(Y)\oplus\widetilde{h}^\bullet(X \wedge Y)\,,
\end{equation}
where $\widetilde{h}$ represents reduced cohomology. We can see that the contribution of $\widetilde{\mathcal{SW}}^5(\bZ_{2^k}[2] \wedge \bZ_{2^k}[2])$ to ${\mathcal{SW}}^5(\bZ_{2^k}[2] \times \bZ_{2^k}[2])$, is nontrivial by comparing it with $\widetilde{\Omega}_{\Spin}^5(\bZ_{2^k}[2] \wedge \bZ_{2^k}[2] )$. Spin cobordism gives the group of maps from the spin bordism groups into $\mathbb{C}^\times$, and when evaluated on a point gives $\Omega^\bullet_{\Spin}(\pt) =\{\bC^\times,\, \bZ_2,\, \bZ_2,\,0,\,\bC^\times,\,0,\,0,\,0,\ldots \} $ in low degrees. 
We claim that it is sufficient to show that $\widetilde{\Omega}_{\Spin}^5(\bZ_{2^k}[2] \wedge \bZ_{2^k}[2] )$ does not vanish. Namely, the bottom three layers of $\Omega^\bullet_{\Spin}(\pt)$ agree with those of $\mathcal{SW}^\bullet(\pt)$, and a fortiori with those of $\SH^\bullet(\pt)$, as was shown in \cite{JF,GJF}. Furthermore, these layers are the only ones that we need to consider in order to compute $\widetilde{\mathcal{SW}}^5(\bZ_{2^k}[2] \wedge \bZ_{2^k}[2])={\SH}^5(\bZ_{2^k}[2] \wedge \bZ_{2^k}[2])$. Here, it is crucial that we are using reduced cohomology. We have a quick-and-dirty way to check that $\widetilde \Omega^5_{\Spin}(\bZ_{2^k}[2] \wedge \bZ_{2^k}[2])$ is nonzero, via the Adams spectral sequence. There are two classes in degree 5, each giving free $\mathcal{A}(1)$-summands in~$\rH^\bullet(B(\bZ_{2^k}[2] \times \bZ_{2^k}[2]); \bZ_2)$, so by Margolis' theorem, the corresponding two $\bZ_2$ summands on the $E_2$-page of the Adams spectral sequence do not admit or receive any differentials. Thus $\widetilde{\Omega}_{\Spin}^5(\bZ_{2^k}[2] \wedge \bZ_{2^k}[2])$ is nontrivial. More generally, this also implies that if $E$ is any group which contains a product of two 2-torsion groups, then the map $\SH^5(E[2])\rightarrow\widetilde{\mathcal{SW}}^5(E[2])$ has non-zero image.

\subsubsection{Bosonic Case} 

In order to classify the bosonic theories associated to braided bosonic strongly fusion 2-categories, it is convenient to work with the associated fermionic theories. This is analogous to how working with an algebra over the real numbers is equivalent to working with the complexified algebra together with the Galois action of $\bZ^T_2$, given by complex conjugation. In this sense, the action of $\bZ^T_2$ provides the necessary data to descend a complex algebra into a real one. The categorification of this classical setup was introduced in \cite{JF2015}. Namely, for symmetric fusion 1-categories, the algebraic closure of $\Vec$ is given by $\SVec$ and the Galois higher group $\Gal(\SVec/\Vec)$ is given by $\bZ^F_2[1]$. This higher group agrees with the physical phenomenon of spin statistics, which says that fermions reverse sign under 360 degree rotation. Then, \textit{Galois descent} asserts that the theory associated to a braided bosonic strongly fusion 2-category $\fB$ is completely described by the $\bZ^F_2[1]$-equivariant theory associated to $\fB\boxtimes 2\SVec$. We can study the later using the equivariant Atiyah-Hirzebruch spectral sequence.

In general, the group of surfaces of $\fB$ is given by a finite abelian group $E$. We begin by showing that $\mathcal{W}^5(\pt)$ does not vanish. That is, we wish to understand the twisted $\mathcal{SW}^\bullet$-cohomology with $E_2$ page given by:
\begin{equation}
    \rH^i(\bZ^F_2[2];\mathcal{SW}^j(\pt))\Rightarrow \mathcal{SW}^{i+j}(\bZ^F_2[2])=\mathcal{W}^{i+j}(\pt).
\end{equation}
To arrive at the last equality, we use the fact that $\mathcal{W}^\bullet$ is the fixed point spectrum of $\mathcal{SW}^\bullet$ under the action of $\bZ^F_2[1]$.
The $E_2$ page is then given by:
\begin{equation}\label{eq:E2braided}
    E^{ij}_2 = \,\begin{array}{c|cccccccccc}
      j\\
      \\
    0&  0  &0&\ldots \\
     \mathcal{SW} &    \mathcal{SW} & 0& \mathcal{SW}_2&\ldots \\
     0& 0&  0 & 0 & 0 & 0 & 0\\
     \bZ_2&\bZ_2&0&\bZ_2&\bZ_2&\bZ_2&\bZ^2_2&\bZ^2_2 &\ldots \\
     \bZ_2&\bZ_2&0&\bZ_2&\bZ_2&\bZ_2&\bZ^2_2&\bZ^2_2 &\ldots \\
     \bC^\times & \bC^\times& 0 & \bZ_2 & 0 & \bZ_4& \bZ_2 & \bZ_2 & \bZ_2 & \bZ_2 \\
     \hline 
     & 0 & 1 & 2 & 3 & 4 & 5 & 6&7&8  & \quad i\,.
\end{array}
\end{equation}
The $d_2$ differentials are the twisted analogue of \eqref{eq:d2differential}
 \begin{align}\label{eq:d2differentialtwisted}
     d_2:& E^{i,2}_2=\rH^i(\bZ_2[2]\,; \bZ_2) \to E^{i+2,1}_2=\rH^{i+2}(\bZ_2[2]\,;\bZ_2) \quad&& X \mapsto \Sq^2 X +\iota_2 X\\\notag 
     d_2:& E^{i,1}_2=\rH^i(\bZ_2[2]\,; \bZ_2) \to E^{i+2,0}_2=\rH^{i+2}(\bZ_2[2]\,; \bC^\times) \quad&& X \mapsto (-1)^{\Sq^2 X+\iota_2 X} \,,
\end{align}
and we find the $E_3$ page is given by 
\begin{equation}
    E^{ij}_3 = \,\begin{array}{c|cccccccccc}
      j\\
      \\
    0&  0  &0&\ldots \\
     \mathcal{SW} &    \mathcal{SW} & 0&  \mathcal{SW}_2&\ldots \\
     0& 0&  0 & 0 & 0 & 0 & 0\\
     \bZ_2&0&0&\bZ_2&0&0&\bZ_2& \ldots \\
     \bZ_2&0&0&0&\bZ_2&0&0& \ldots \\
     \bC^\times & \bC^\times& 0 & 0 & 0 & \bZ_4& \bZ_2 & 0 & 0 & 0 \\
     \hline 
     & 0 & 1 & 2 & 3 & 4 & 5 & 6&7&8  & \quad i\,.
\end{array}
\end{equation}
The $d_5$ differential from $(0,4)$ records the obstruction to minimal modular extensions. It sends a class in $\mathcal{SW}$ to 0 if the minimal modular extension exists, and to 1 if it does not exist \cite{JF}.
The main result of \cite{JFR} shows that the possible $d_5$ vanishes.
We therefore find that $\mathcal{W}^{5}(\pt)\cong \bZ_2$.

Now, let us consider any finite abelian group $E$. It follows from \eqref{eq:prodgroup}, that $$\mathcal{W}^{5}(E[2])\cong \mathcal{W}^{5}(\pt) \oplus \widetilde{\mathcal{SW}}^{5}(E[2])\oplus \widetilde{\mathcal{SW}^5}(\bZ^F_2[2]\wedge E[2]).$$ It follows from what we have argued above in the fermionic case that the canonical map~$\rH^5(E[2]; \mathbb{C}^{\times})\rightarrow \widetilde{\mathcal{SW}}^{5}(E[2])$ is non-zero for a general finite abelian group $E$. As a consequence, the theory associated to a braided bosonic strongly fusion 2-category can not be condensed to the vacuum in general.

Before moving on the case of symmetric strongly fusion 2-categories, let us briefly remark that, in section \ref{section:ExampleCats}, we have also considered examples when the condensation yields a sylleptic strongly fusion 2-category.
The computations for the theories associated to these 2-categories were performed in \cite{JFY2}, where the object of study was topological (4+1)d theories.

\subsection{Symmetric Strongly Fusion 2-Categories}

We now analyze the structure of symmetric strongly fusion 2-categories. More precisely, we will show below that every symmetric strongly fusion 2-category admits a fibre 2-functor to $2\SVec$. In the process, we will also show that every symmetric fermionic strongly fusion 2-category is completely determined by its groups of connected components. These computations establish the 2-Deligne theorem for symmetric fusion 2-categories. Namely, it follows from corollary \ref{cor:stronglysylleptic} together with the obvious fermionic analogue, that every symmetric fusion 2-category admits a fibre 2-functor to a strongly fusion 2-category. Putting the above discussion together, we obtain the following theorem, which is a categorification of \cite{deligne2002}.

\begin{theorem}\label{thm:2Deligne}
Every symmetric fusion 2-category admits a fibre 2-functor to $2\SVec$.
\end{theorem}

\noindent We point out that this result was first announced in \cite{JF3}. In addition, we expect that the above theorem can be used to classify symmetric fusion 2-categories. More precisely, every symmetric fusion 2-category should be equivalent to the symmetric monoidal 2-category of finite semisimple 2-representation of a ``super 2-group''.


\subsubsection{Fermionic Case}\label{subsub:fermionicsymmetric}

Let $\mathfrak{S}$ be a symmetric fermionic strongly fusion 2-category, and let us denote by $E$ its abelian group of connected components. We now wish to understand what additional data besides $E$, if any, needs to be supplied to recover $\mathfrak{S}$. We begin by describing $\mathfrak{S}^{\times}$ the Picard sub-2-category of $\mathfrak{S}$, that is the maximal sub-2-category on the invertible objects and morphisms.

It has been established in \cite{GJO} that the homotopy theory of symmetric monoidal 2-categories for which all objects and morphisms are invertible is equivalent to that of spectra with homotopy groups concentrated in degrees $0,1,$ and $2$. In particular, the Picard 2-category $\mathfrak{S}^{\times}$ fits into the following fibre sequence of spectra $$2\SVec^{\times}\rightarrow \mathfrak{S}^{\times}\rightarrow  \mathrm{H}E\rightarrow \Sigma 2\SVec^{\times},$$ where $\mathrm{H}E$ denotes the Eilenberg-MacLane spectrum associated to $E$. In particular, $\mathfrak{S}^{\times}$ is completely determined by the map of spectra $\mathrm{H}E\rightarrow \Sigma 2\SVec^{\times}$. Up to homotopy, such maps are classified by the group $\SH^{7}(E[4])$.

In order to compute the group $\SH^{7}(E[4])$, we invoke the Atiyah-Hirzebruch spectral sequence with the $E_2$-page:
\begin{equation}\label{eq:SpinBord}
    \rH^i(E[4];\,\SH^j(\pt))\Longrightarrow \SH^{i+j}(E[4]) \,.
\end{equation}
We will show that the degree seven supercohomology group $\SH^{7}(E[4])$ vanishes for any finite abelian group $E$. Firstly, it follows from \cite{EilenbergMacLane} that $\rH^7(E[4]; \bC^\times)=0$ if $E$ has no 2-torsion. In addition, the Hurewicz theorem shows that $\SH^{7}(E[4])$ can only be non-trivial if $E$ has 2-torsion. We start with the case $E=\bZ_{2^k}$, as explained in \cite{Serre}, the cohomology $\rH^{\bullet}(\bZ_{2^k}[n],\bZ_2)$ is a polynomial ring $\bZ_2[\Sq^I(\iota_n)]$
where the generator $\iota_n \in \rH^n(\bZ_{2^k}[n]; \bZ_2)$ is in degree $n$, and $I=(i_1,i_2,\ldots,i_m)$ runs over all sequences such that $i_j \geq 2i_{j+1}$
of excess $e(I)<n$. This quantity is defined as $e(I) = i_2 - \sum_{j \geq 3 } i_j$ and  $\Sq^I x = \Sq^{i_1} \Sq^{i_2} \ldots \Sq^{i_m} x$. If $i_m =1$ then $\Sq^I x =\Sq^{i_1}\Sq^{i_2}\ldots \Sq^{i_{m-1}} \beta_k x $ where $\beta_k$ denotes the $k$-th power Bockstein for the short exact sequence
\begin{equation}
    0 \rightarrow \bZ_2 \rightarrow \bZ_{2^{k+1}} \rightarrow \bZ_{2^k} \rightarrow 0\,.
\end{equation}
The $E_2$ page for \eqref{eq:SpinBord} in terms of the generators then takes the form
\begin{equation}
    E^{ij}_2 = \,\begin{array}{c|ccccccccc}
      j\\
      \\
     \bZ_2&\bZ_2&0&0&0&\iota_4 &\beta_k\iota_4  & \Sq^2  \iota_4& \ldots \\
     \bZ_2&\bZ_2&0&0&0&\iota_4 &\beta_k\iota_4 & \Sq^2 \iota_4& (\Sq^3 \iota_4, \,\Sq^2 \beta_k \iota_4) \\
     \bC^\times & \bC^\times& 0 & 0 & 0 & (-1)^{\iota_4}& 0 &(-1)^{ \Sq^2 \iota_4} & (-1)^{\Sq^2 \beta_k \iota_4}  \\
     \hline 
     & 0 & 1 & 2 & 3 & 4 & 5 &  6& 7&8\quad i\,.
\end{array}
\end{equation}
The $d_2$ differentials are given by:
 \begin{align}\label{eq:d2differential}
     d_2:& E^{i,2}_2=\rH^i(\bZ_{2^k}[4]\,; \bZ_2) \to E^{i+2,1}_2=\rH^{i+2}(\bZ_{2^k}[4]\,;\bZ_2) \quad&& X \mapsto \Sq^2 X \\\notag 
     d_2:& E^{i,1}_2=\rH^i(\bZ_{2^k}[4]\,; \bZ_2) \to E^{i+2,0}_2=\rH^{i+2}(\bZ_{2^k}[4]\,; \bC^\times) \quad&& X \mapsto (-1)^{\Sq^2 X} \,,
\end{align}
and there are $d_2$'s leaving the entries in bidegrees $(4,2)$ and $(4,1)$ that carry the generator $\iota_4$ to $\Sq^2 \iota_4$ and are isomorphisms. Additionally, there are $d_2$ differentials leaving the entries in bidegrees $(5,1)$ and $(5,2)$ which are isomorphisms.   
In total degree seven, the $E_3$ page converges to the $E_\infty$ page and we see that $\SH^7(\bZ_{2^k}[4])=0$. 
If $E$ is a product of groups, we can use \eqref{eq:prodgroup}, where the spaces are fourth Eilenberg-MacLane spaces of groups that are 2-torsion. Then the term corresponding to $\widetilde{h}^\bullet(X \wedge Y)$ for supercohomology will only begin to contribute in degree 8, and everything else vanishes. In summary, we have shown that $\SH^{7}(E[4])=0$ for any group $E$.

This implies that $\mathfrak{S}^{\times} \cong 2\SVec^{\times} \times E$ as symmetric monoidal 2-categories. In particular, $\mathrm{B}\SVec \times E$ is a full symmetric monoidal sub-2-category of $\mathfrak{S}$. But $\mathrm{B}\SVec \times E$ contains an object in every connected component of $\mathfrak{S}$, so that its Cauchy completion $Cau(\mathrm{B}\SVec \times E)\simeq 2\SVec[E]$ is equivalent to $\mathfrak{S}$ as a symmetric monoidal 2-category. Thus, we obtain the following result.

\begin{proposition}
Every symmetric strongly fusion 2-category is of the form $2\SVec[E]$ for some finite abelian group $E$.
\end{proposition}

\noindent In particular, every symmetric strongly fusion 2-category admits a fibre 2-functor to $2\SVec$.


\subsubsection{Bosonic Case}

For a symmetric bosonic strongly fusion 2-category, the obstruction to condensing to the symmetric fusion 2-category to $2\Vec$ is given by a class in $\rH^7(E[4];\bC^\times)$.
The group $\rH^{n+m+1}(E[n]; \bC^\times)$ may be thought of as parameterizing the ways for $m$ spacetime dimensional objects to fuse in $n$ ambient dimensions with fusion rule $E$.
A computation of this cohomology group can be found in \cite{EilenbergMacLane} where the authors obtained $\rH^7(E[4];\bC^\times) = \widehat{({E}_2)}$\,, with $E_2$ the 2-torsion subgroup of $E$, and for a group $A$ we denote $\widehat{A} = \hom(A, \mathbb{C}^\times)$. Even though this cohomology group is not necessarily trivial, the computation in \ref{subsub:fermionicsymmetric} shows that if we work in a fermionic setting, then there is no obstruction to the existence of a fibre 2-functor to $2\SVec$. In particular, any symmetric bosonic strongly fusion 2-category admits a fibre 2-functor to $2\SVec$, which concludes the proof of theorem \ref{thm:2Deligne}.

\appendix

\section{Diagrams}

\subsection{Proof of Theorem \ref{thm:SymmetricAlgebraModule}}\label{subsub:proofSymmetricAlgebraModule}

\begin{figure}[!b]
\centering
\includegraphics[width=160mm]{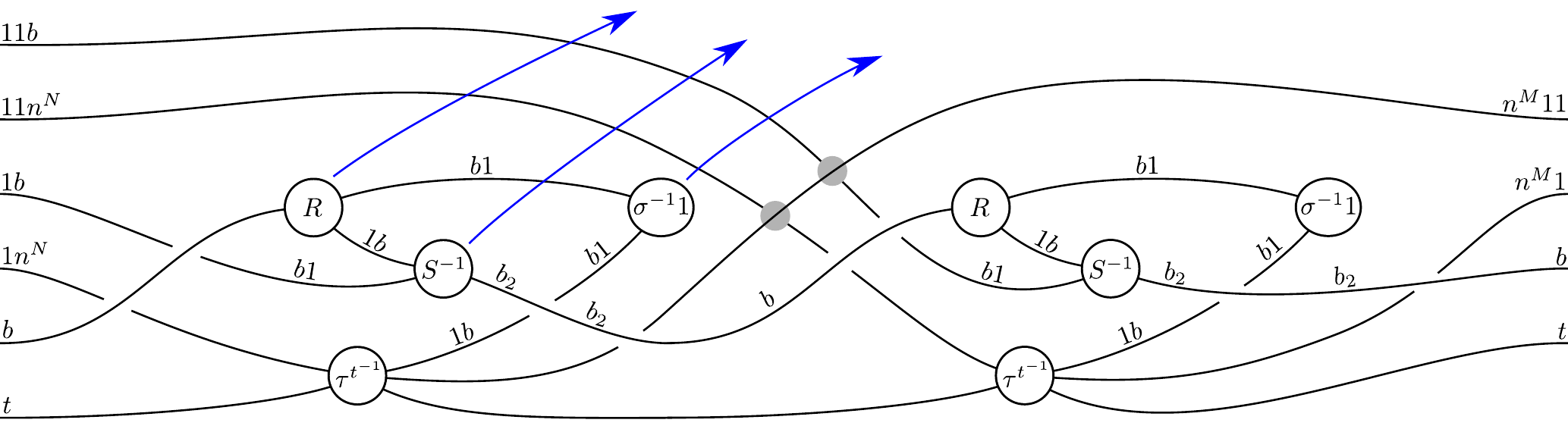}
\caption{Axiom a (Part 1)}
\label{fig:taubbalanced1}
\end{figure}

\begin{figure}[!hbt]
\centering
\includegraphics[width=160mm]{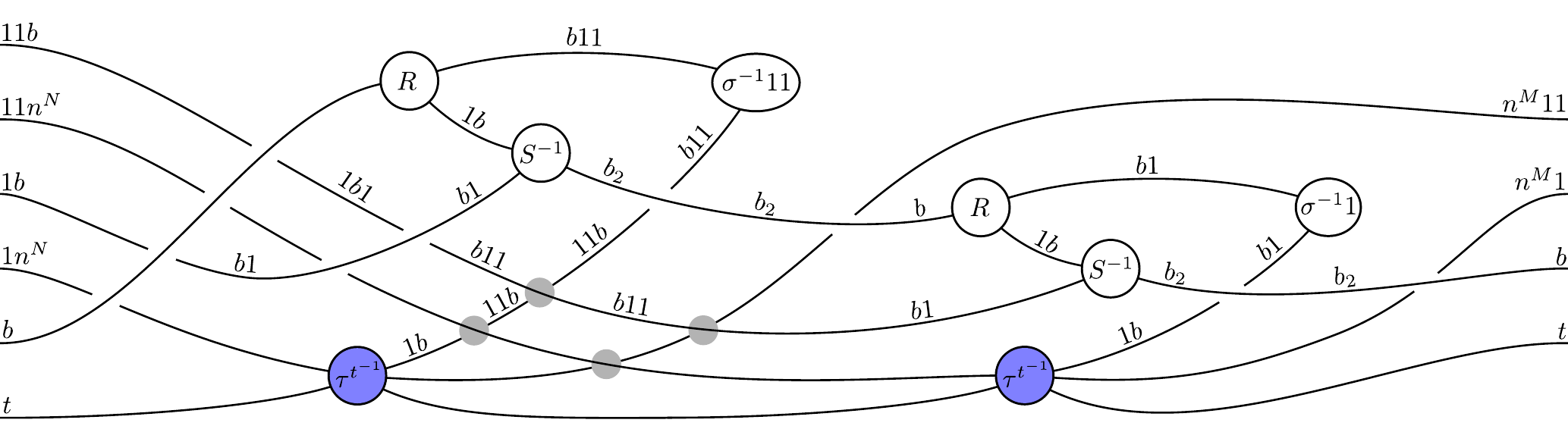}
\caption{Axiom a (Part 2)}
\label{fig:taubbalanced2}
\end{figure}

\begin{figure}[!hbt]
\centering
\includegraphics[width=160mm]{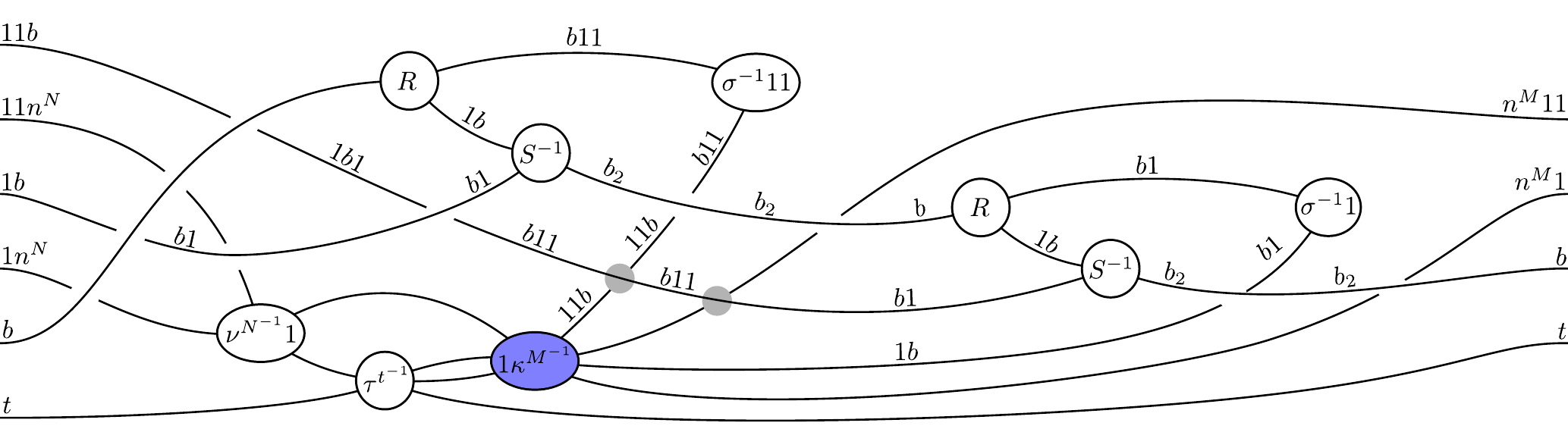}
\caption{Axiom a (Part 3)}
\label{fig:taubbalanced3}
\end{figure}

\begin{figure}[!hbt]
\centering
\includegraphics[width=165mm]{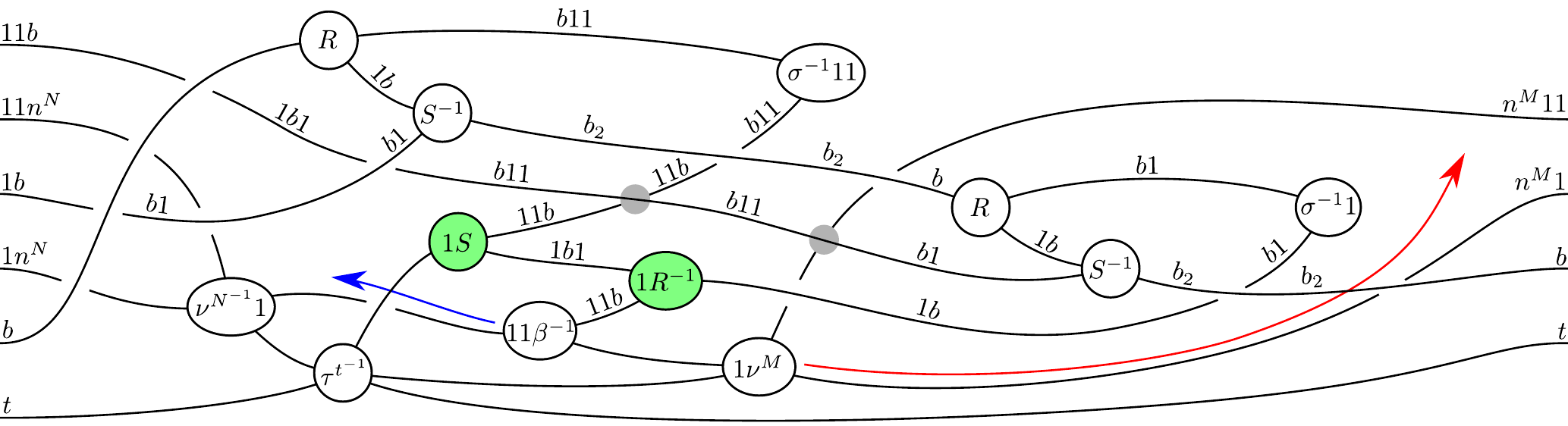}
\caption{Axiom a (Part 4)}
\label{fig:taubbalanced4}
\end{figure}

\begin{figure}[!hbt]
\centering
\includegraphics[width=165mm]{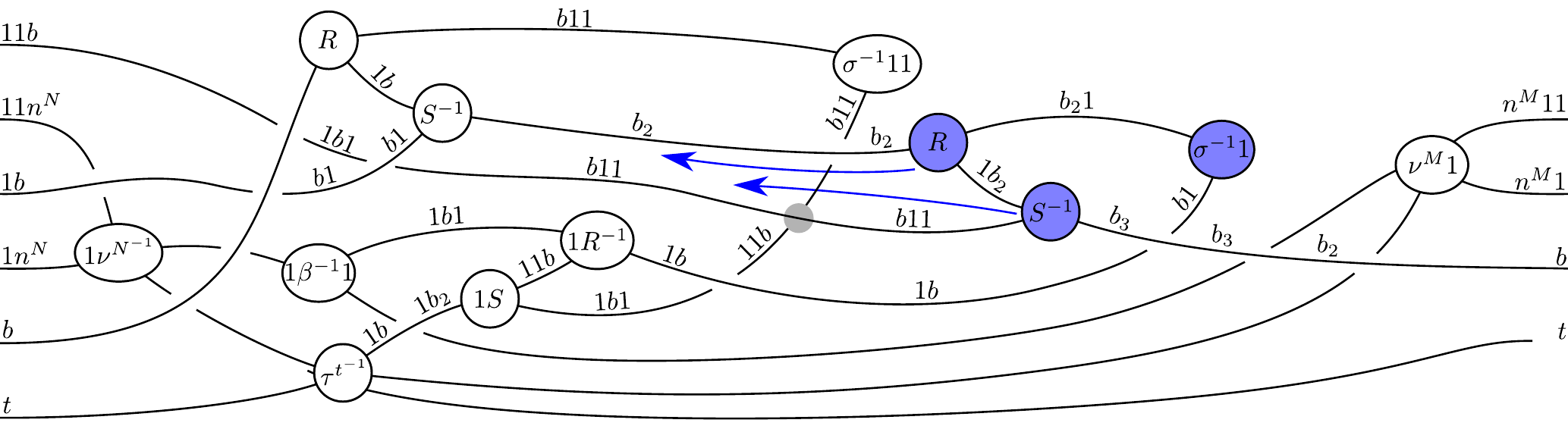}
\caption{Axiom a (Part 5)}
\label{fig:taubbalanced5}
\end{figure}

\begin{figure}[!hbt]
\centering
\includegraphics[width=165mm]{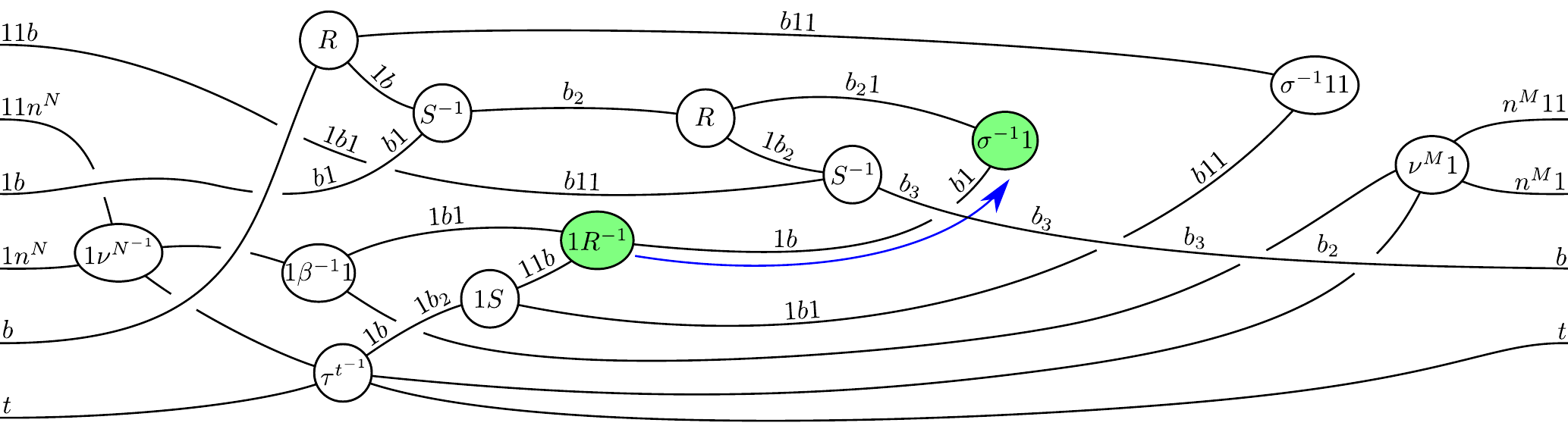}
\caption{Axiom a (Part 6)}
\label{fig:taubbalanced6}
\end{figure}

\begin{figure}[!hbt]
\centering
\includegraphics[width=165mm]{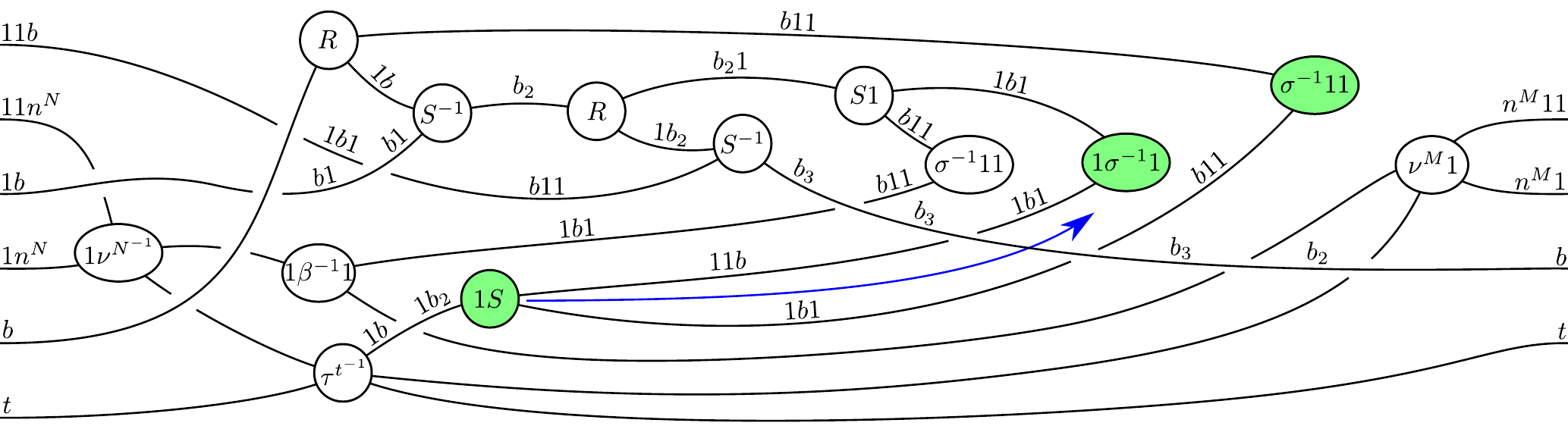}
\caption{Axiom a (Part 7)}
\label{fig:taubbalanced7}
\end{figure}

\begin{figure}[!hbt]
\centering
\includegraphics[width=165mm]{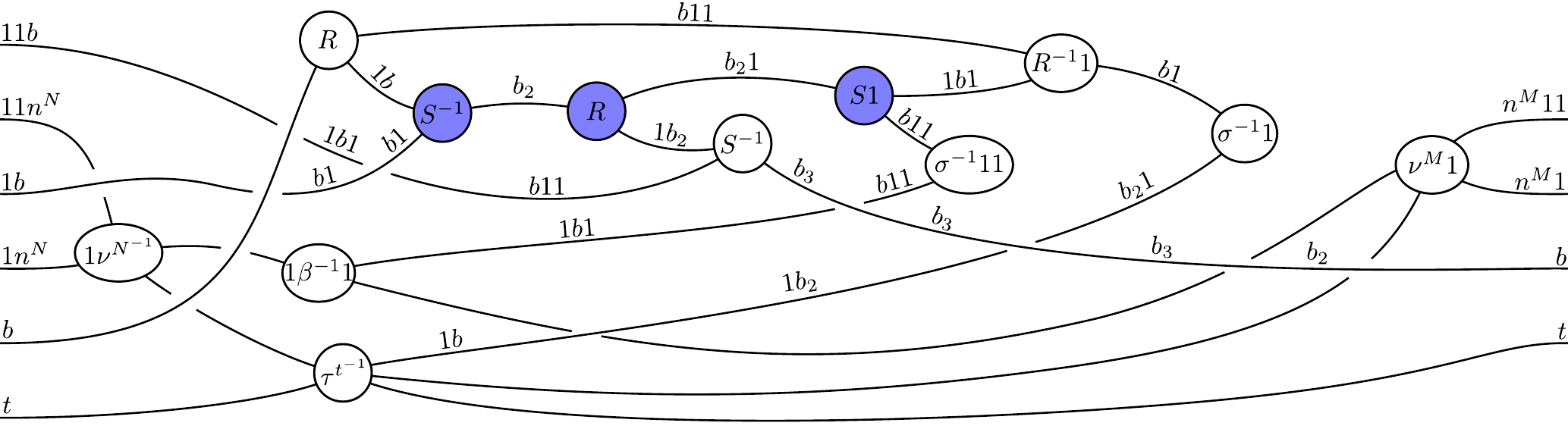}
\caption{Axiom a (Part 8)}
\label{fig:taubbalanced8}
\end{figure}

\begin{figure}[!hbt]
\centering
\includegraphics[width=165mm]{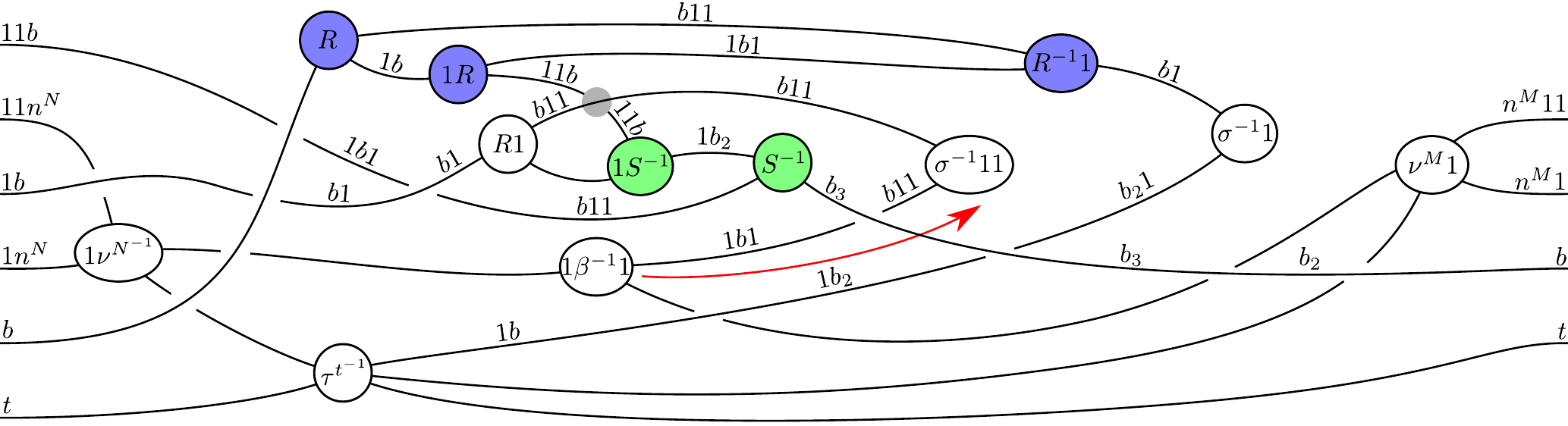}
\caption{Axiom a (Part 9)}
\label{fig:taubbalanced9}
\end{figure}

\begin{figure}[!hbt]
\centering
\includegraphics[width=165mm]{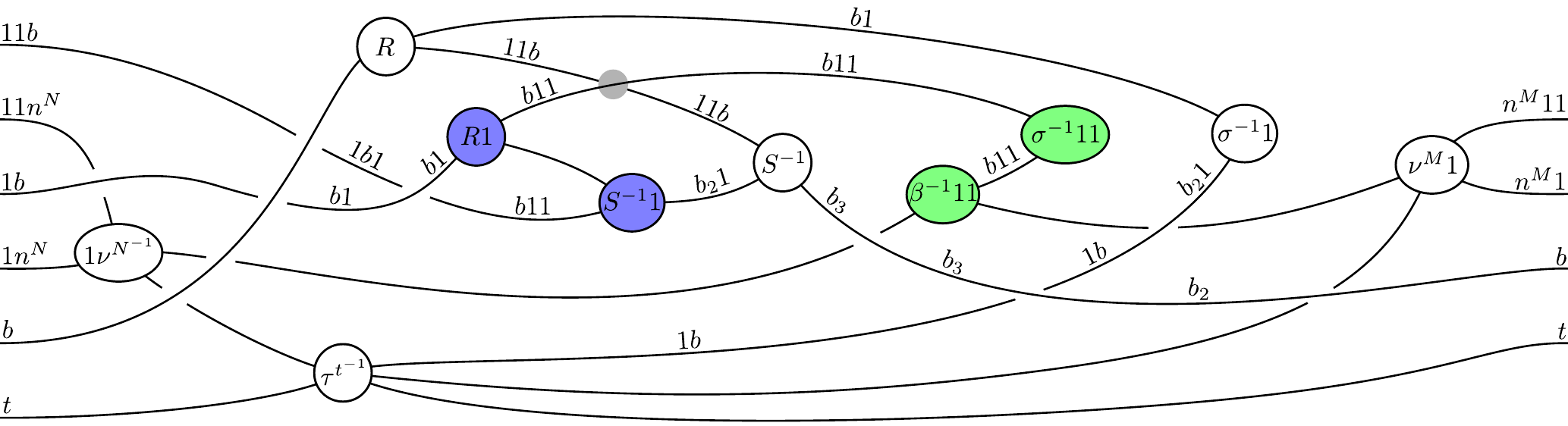}
\caption{Axiom a (Part 10)}
\label{fig:taubbalanced10}
\end{figure}

\begin{figure}[!hbt]
\centering
\includegraphics[width=165mm]{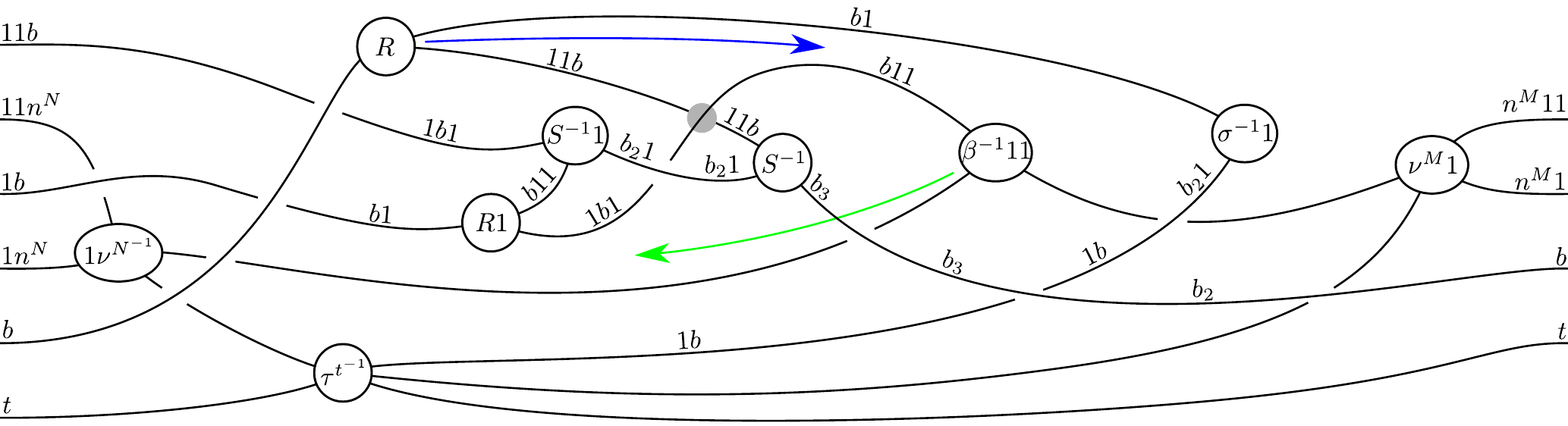}
\caption{Axiom a (Part 11)}
\label{fig:taubbalanced11}
\end{figure}

\begin{figure}[!hbt]
\centering
\includegraphics[width=120mm]{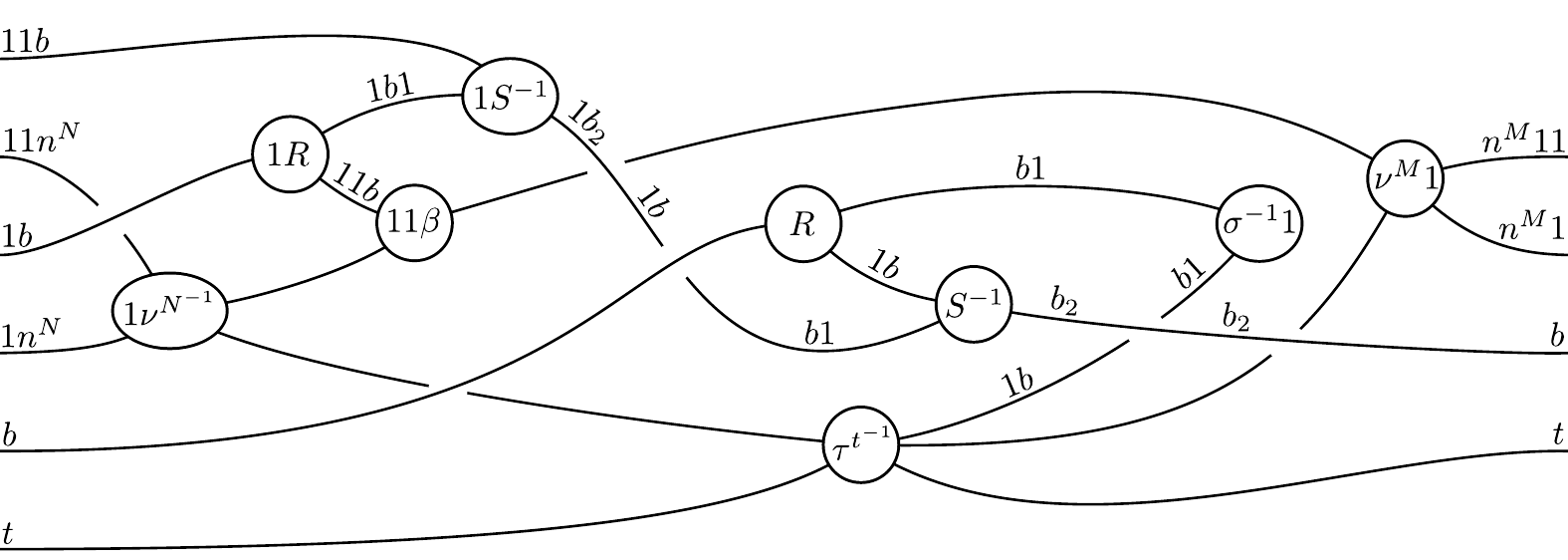}
\caption{Axiom a (Part 12)}
\label{fig:taubbalanced12}
\end{figure}

\FloatBarrier

\subsection{Proof of Proposition \ref{prop:SymmetricAlgebraCenterModule}}\label{subsub:proofSymmetricAlgebraCenterModule}

\begin{figure}[!b]
\centering
\begin{minipage}{.5\textwidth}
  \centering
    \includegraphics[width=45mm]{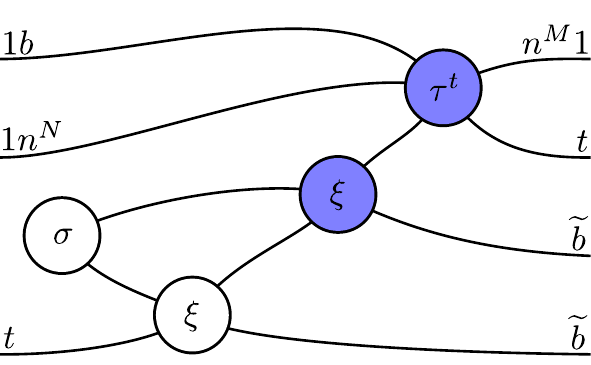}
    \caption{Balanced structure (Part 1)}
    \label{fig:sigmaxixibalanced1}
\end{minipage}%
\begin{minipage}{.5\textwidth}
  \centering
    \includegraphics[width=45mm]{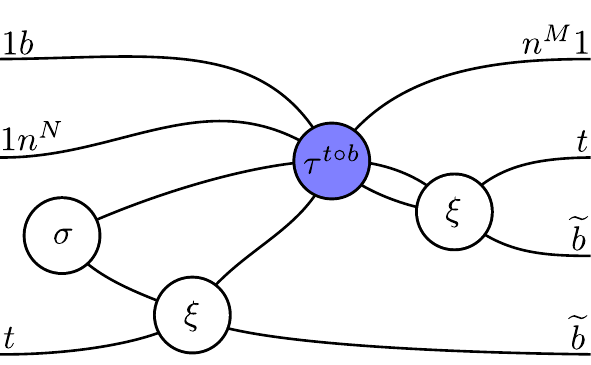}
    \caption{Balanced structure (Part 2)}
    \label{fig:sigmaxixibalanced2}
\end{minipage}
\end{figure}

\begin{figure}[!hbt]
\centering
\begin{minipage}{.5\textwidth}
  \centering
    \includegraphics[width=75mm]{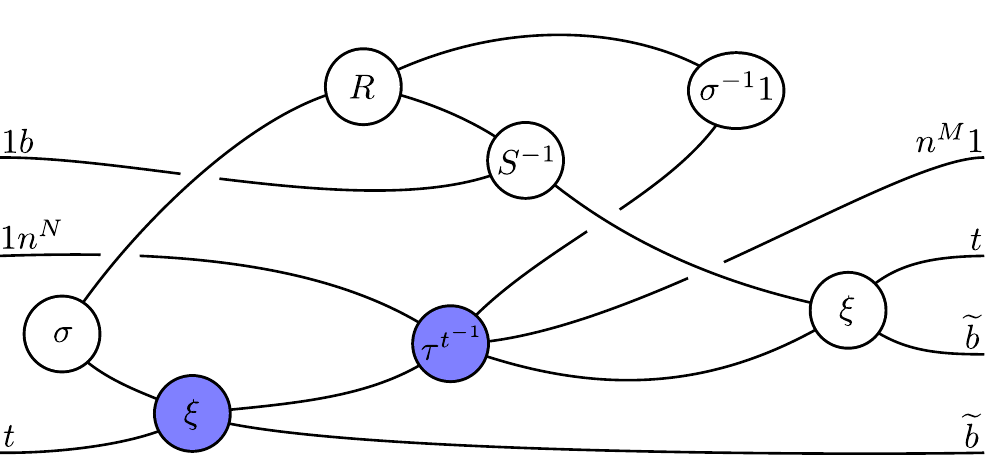}
    \caption{Balanced structure (Part 3)}
    \label{fig:sigmaxixibalanced3}
\end{minipage}%
\begin{minipage}{.5\textwidth}
  \centering
    \includegraphics[width=75mm]{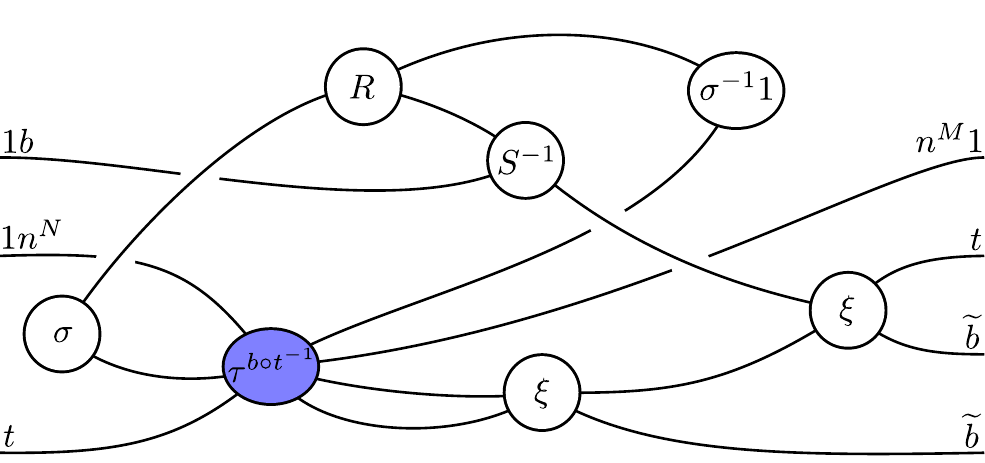}
    \caption{Balanced structure (Part 4)}
    \label{fig:sigmaxixibalanced4}
\end{minipage}
\end{figure}

\begin{figure}[!hbt]
\centering
\includegraphics[width=105mm]{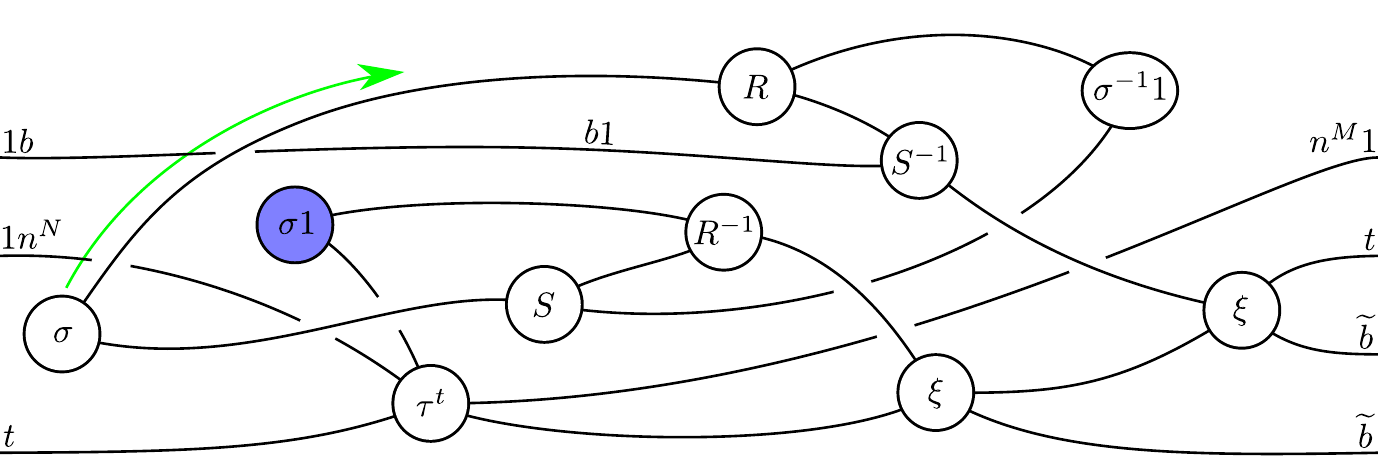}
\caption{Balanced structure (Part 5)}
\label{fig:sigmaxixibalanced5}
\end{figure}

\begin{figure}[!hbt]
\centering
\includegraphics[width=90mm]{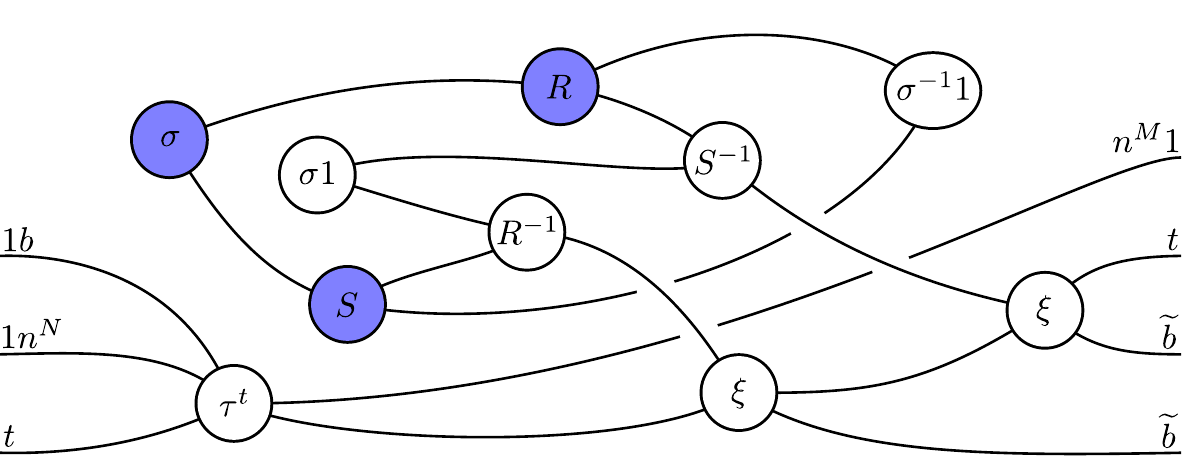}
\caption{Balanced structure (Part 6)}
\label{fig:sigmaxixibalanced6}
\end{figure}

\begin{figure}[!hbt]
\centering
\begin{minipage}{.67\textwidth}
  \centering
    \includegraphics[width=90mm]{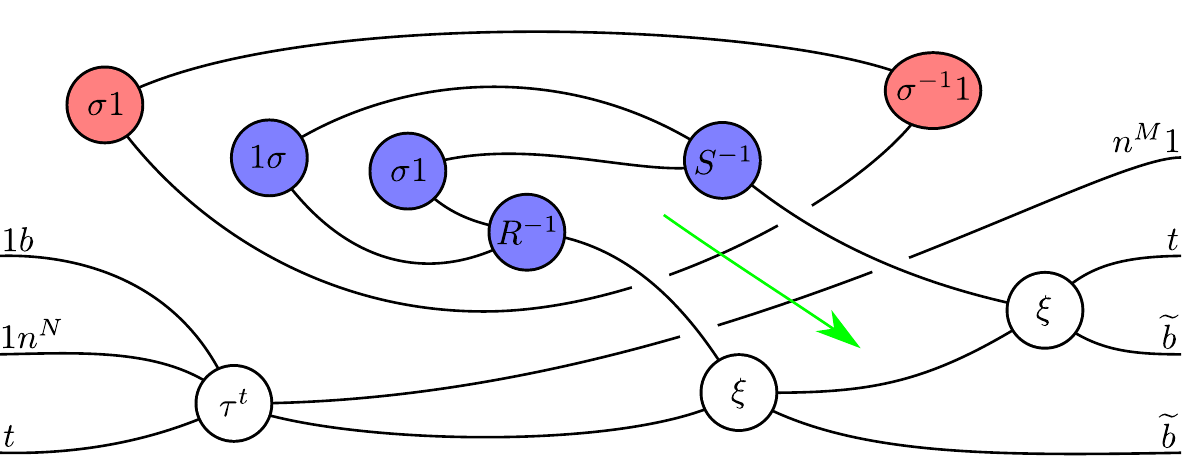}
    \caption{Balanced structure (Part 7)}
    \label{fig:sigmaxixibalanced7}
\end{minipage}%
\begin{minipage}{.33\textwidth}
  \centering
    \includegraphics[width=45mm]{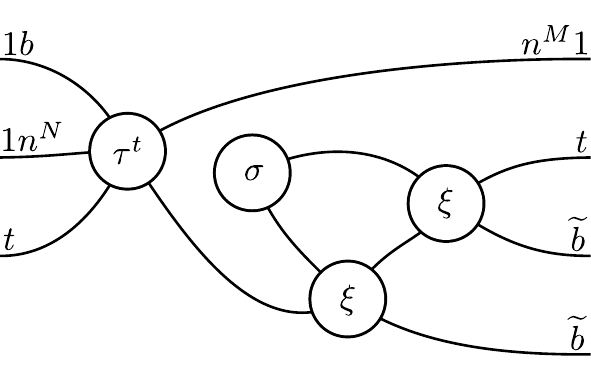}
    \caption{Balanced structure (Part 8)}
    \label{fig:sigmaxixibalanced8}
\end{minipage}
\end{figure}

\FloatBarrier

\bibliographystyle{alpha}
\bibliography{references}{}
\end{document}